\documentclass{gtpart}  
\usepackage{amscd}
\usepackage{amsgen}
\usepackage{curvesls}
\usepackage{epic}
\usepackage[all]{xy}
\usepackage[dvips]{graphicx}
\usepackage{psfrag}

\theoremstyle{plain}
\newtheorem{prop}[subsection]{Proposition}
\newtheorem{thm}[subsection]{Theorem}
\newtheorem{lem}[subsection]{Lemma}
\newtheorem{cor}[subsection]{Corollary}
\newtheorem{eg}[subsection]{Example}

\newtheorem{defn}[subsection]{Definition}

\theoremstyle{remark}

\theoremstyle{definition}
\newtheorem{exm}[subsection]{Example}

\numberwithin{equation}{section}

\newcommand{\Real}{{\mathbb R}}

\newcommand{\K}{{\mathcal{K}}}

\newcommand{\Rational}{{\mathbb Q}}
\newcommand{\Cu}{{\mathcal{C}}}
\newcommand{\Diff}{{\mathrm{Diff}}}
\newcommand{\Emb}{{\mathrm{Emb}}}

\newcommand{\EK}[1]{{\mathrm{EC}({#1})}}
\newcommand{\cyc}[1]{{\Zed/#1\Zed}}
\newcommand{\EC}[2]{{Emb(\Real^{#1} \times {#2},\Real^{#1} \times {#2})}}

\newcommand{\splice}{{\bowtie}}
\newcommand{\rep}{{\Delta}}

\newcommand{\I}{{\mathsf I}}
\newcommand{\Seifert}{{\mathcal S}}
\newcommand{\Keychain}{{\mathcal H}}
\newcommand{\IG}{{\Bbb G}}

\newcommand{\cfg}{{\mathrm{Conf}}}

\newcommand{\Natural}{{\mathbb N}}

\newcommand{\Zed}{{\mathbb Z}} 
\newcommand{\Complex}{{\mathbb C}}
\newcommand{\F}{{\mathbb F}}

\begin{document}

\title[On the homology of the space of knots] {On the homology of the space of knots}
\author{R.~Budney}
\address{
Department of Mathematics and Statistics \\
University of Victoria \\
Victoria BC Canada \\
V8W 3P4} \email{rybu@uvic.ca}
\author{F. R.~Cohen*}
\address{Department of Mathematics \\
University of Rochester\\ Rochester, NY 14627 U.S.A.}
\email{cohf@math.rochester.edu}
\thanks{*Partially supported by the NSF Grant No. DMS-0072173 and CNRS-NSF
Grant No. 17149}

\subject{Primary: 58D10, 57T25; Secondary: 57M25,57Q45}

\begin{abstract}
Consider the space of `long knots' in $\Real^n$ $\K_{n,1}$. This is
the space of knots as studied by V.~Vassiliev. Based on previous
work \cite{B,CLM}, it follows that the rational homology of
$\K_{3,1}$ is free Gerstenhaber-Poisson algebra. A
partial description of a basis is given here. In addition, the
mod-$p$ homology of this space is a `free, restricted
Gerstenhaber-Poisson algebra'. Recursive application of this theorem
allows us to deduce that there is $p$-torsion of all orders in the
integral homology of $\K_{3,1}$.

This leads to some natural questions about the homotopy type of the
space of long knots in $\Real^n$ for $n>3$, as well as consequences
for the space of smooth embeddings of $S^1$ in $S^3$ and embeddings of
$S^1$ in $\Real^3$.
\end{abstract}

\maketitle

\section{Introduction}\label{INTRODUCTION}

The purpose of this paper is to give homological properties of the
classical spaces of smooth `long' embeddings $\K_{3,1} =
\Emb(\Real,\Real^3) $ and smooth embeddings $\Emb(S^1,S^3)$. Some
results here also apply to the embedding spaces $\Emb(S^j,S^n)$, and
`long' embedding spaces $\K_{n,j}= \Emb(\Real^j,\Real^n)$, with the
main results focused on the $3$-dimensional case $j=1$, $n=3$.

The approach here to these homological problems follows recent work
of Hatcher \cite{Hatcher3, Hatcher4} and Budney \cite{B, B2}. The
homotopy type of the components of $\Emb(S^1,S^3)$ and $\K_{3,1}$
are understood completely in terms of configuration spaces in the
plane, Stiefel manifolds, isometry groups of certain hyperbolic link
complements and various natural iterated bundle operations. Many of
the homological properties of both $\K_{3,1}$, and $\Emb(S^1,S^3)$
follow from combining this information together with earlier work of
Cohen \cite{CLM} on configuration spaces.

The space $\K_{n,1}$ admits the structure of an $H$-space induced by
concatenation of `long' embeddings. In addition, this $H$-space
structure for $\K_{3,1}$ was shown to extend to a free $\Cu_2$-space
in the sense of May, with ``generating set" given by the space of
prime long knots \cite{B}.

One consequence is that the homotopy type of the space of long knots
is determined completely by the homotopy type of the prime long
knots. Information concerning spaces of prime long knots is combined
with bundle theoretic constructions to give a large contribution to
the homology groups for spaces of long knots, as well as
$\Emb(S^1,S^3)$.

The structure of a graded analogue of a Poisson algebra, a
Poisson-Gerstenhaber algebra, arises in the work here. An
introduction to Poisson algebras is given in \cite{CP}, pages
$177$-$182$ while some applications are given in \cite{CLM} pages
$215$-$216$ and \cite{Tourtchine}. A Poisson-Gerstenhaber algebra
$A$ is a graded commutative algebra over $\mathbb Q$ given by $A_n$
in gradation $n$ together with a graded skew symmetric bilinear map
$$\{-,-\}:A_s \otimes A_t \to A_{s+t+1}$$ which satisfies the
following where $|a|$ denotes the degree of an element $a$ in $A$:

\begin{enumerate}
    \item the Jacobi identity
$$\{a,\{b,c\}\} = (-1)^{|a|+|b|+|c|+1}\{\{a,b\},c\} + (-1)^{|b||c| + |c|+1}\{\{a,c\},b\}$$
where the signs will be typically omitted with the above written as
$$\{\{a,\{b,c\}\} = (\pm 1)\{\{a,b\},c\} + (\pm 1)\{\{a,c\},b\}.$$
    \item the Leibniz formula
$$\{a\cdot b,c\} = a\cdot \{b,c\} + (-1)^{|b||a|}b\{a,c\}.$$
\end{enumerate}

A standard example of such a Gerstenhaber-Poisson algebra is given
by the rational homology algebra of $\Omega^2(X)$ for $X$ a bouquet
of spheres of dimension at least $3$ by \cite{CLM} in which the
precise axioms are recorded.

The mod-$p$ homology of $\K_{3,1}$ has more detailed structure, and
is, loosely speaking, a `free restricted Gerstenhaber-Poisson
algebra' with additional structure satisfied by free $\mathcal
\Cu_2$-spaces \cite{CLM}, freely generated on the mod-$p$ homology
of the subspace of prime long knots.

The main results of the current article are Theorem \ref{homology of
long knots in 3 space} on the structure of the homology of
$\K_{3,1}$, Proposition \ref{prop:embeddings of circles in the three
sphere} concerning implications for the space of smooth embeddings
of $S^1$ in $S^3$, Proposition \ref{prop:Two torsion in h1}, a
homological characterization of the unknot, as well as Theorem
\ref{first_oddp_torsion} on the minimal $i$ such that
$H_i(\K_{3,1};\Zed)$ contains $p$-torsion.

\begin{defn}
$ \K_{n,j} := \{ f : \Real^j \to \Real^n : f \text{ is an embedding
and } f(x_1,x_2,\cdots,x_j)=(x_1,x_2,\cdots,x_j,0,\cdots, 0) \text{
for } |x| \geq 1 \}$. $\K_{n,1}$ is traditionally called the space
of long knots in $\Real^n$, and $\K_{n,j}$ the space of long
$j$-knots in $\Real^n$. Given an element $f \in \K_{n,j}$ the
connected-component of $\K_{n,j}$ containing $f$ is denoted
$\K_{n,j}(f)$. Two knots are considered equivalent if they are in
the same connected component of $\K_{n,j}$ (That is the knots are
isotopic).
\end{defn}

Let $X_{\K} = \{ f \in \K_{3,1} : f \text{ is prime} \}$, where the
word `prime' is used in the traditional sense of Schubert
\cite{Sch}, that is $X_{\K}$ is the union of the connected
components of $\K_{3,1}$ which contain knots that are not connected
sums of two or more non-trivial knots, nor are they allowed to
contain the unknot.

\begin{thm}\label{long knots in 3 space} \cite{B1}
The space $\K_{3,1}$ is homotopy equivalent to
$$C(\Real^2,X_{\K} \amalg \{*\} )$$
that is the labelled configuration space of points in the plane with
labels in $X_{\K} \amalg \{*\} $. Furthermore, the following hold:
\begin{enumerate}
    \item Each path-component of
$\K_{3,1}$ is a $K(\pi,1)$.
    \item The path-components of $\Cu_2(n) \times_{\Sigma_n}(X_{\K})^n$
    for all $n$, and thus the
path-components of $\K_{3,1}$ are given by
$$\Cu_2(n) \times_{\Sigma_f} \prod_{i=1}^n \K_{3,1}({f_i}) $$
for certain choices of $f_1,\cdots, f_n \in X_{\K}$, and Young
subgroups $\Sigma_f$.
\end{enumerate}
\end{thm}

The above theorem can be thought of as a generalization of
Schubert's Theorem which states that $\pi_0 \K_{3,1}$ is a free
commutative monoid on countably-infinite many generators \cite{Sch}.
Schubert's theorem is about the monoid structure on $\pi_0 \K_{3,1}$
induced by the cubes action, while the above theorem is space-level
on $\K_{3,1}$.

In general, $\K_{n,1}$ is a homotopy-associative $H$-space with
multiplication induced by concatenation. This multiplication gives a
product operation $$H_s(\K_{n,1}) \otimes H_t(\K_{n,1}) \to
H_{s+t}(\K_{n,1})$$ Since $\K_{3,1}$ admits the action of the operad
of little $2$-cubes, there is an induced map
$$\theta: S^1 \times \K_{3,1} \times \K_{3,1} \to \K_{3,1}$$
together with an operation in homology with any coefficients
$$H_s(\K_{3,1}) \otimes H_t(\K_{3,1}) \to H_{1+s+t}(\K_{3,1})$$
which is denoted, up to sign, by
$$\{\alpha,\beta\} \equiv \lambda_1(\alpha,\beta) = \theta_*(\iota \otimes \alpha \otimes \beta)$$
for $\alpha$ in $H_s(\K_{3,1})$, $\beta$ in $H_t(\K_{3,1})$ and
$\iota \in H_1(S^1)$ the fundamental class. These operations satisfy
the structure of a graded Poisson algebra for which the bracket
operation $\lambda_1(\alpha,\beta)$ is called the Browder operation
in \cite{CLM}.

The next result uses the product operation above as well as the
bracket operation $\{\alpha,\beta\} = \lambda_1(\alpha,\beta)$, and
follows by interweaving the results of \ref{long knots in 3 space},
and \cite{CLM}.  We will use the notation $\F_p = \Zed /p\Zed$
to denote the field with $p$ elements, when $p$ is a prime. 
To state these results, additional information
given by three functors from graded vector spaces $V$ over a field
$\F$ are described next with complete details given in Section
\ref{c2homol}:

\begin{enumerate}
    \item If the characteristic of the field is $0$ then the value of
    the functor on objects $V$ is the symmetric algebra generated
    by an algebraically `desuspended' free Lie algebra generated
    by the suspension of $V$, and denoted
    $$S[\sigma^{-1}L[\sigma(V)]].$$ This last algebra is a free
    Gerstenhaber-Poisson algebra.

    \item For the field $\F_2$ the value of
    the functor on objects $V$ is the symmetric algebra generated
    by an algebraically `desuspended'  free, mod-$2$ restricted Lie algebra
    generated by the suspension of $V$, and denoted
    $$S[\sigma^{-1}L^{(2)}[\sigma(V)]].$$ This last algebra is a
    graded version of a free restricted Lie algebra.

    \item For the field $\F_p$ ($p$ an odd prime), then the value of
    the functor on objects $V$ is the symmetric algebra generated
    by an algebraically ``desuspended"  mod-$p$ free restricted Lie algebra
    generated by the suspension of $V$ plus an additional summand as described
    in Section \ref{c2homol}, and denoted
    $$S[\sigma^{-1}L^{(p)}[\sigma(V)] \oplus \sigma^{-2}W^p[\sigma(V)]].$$
\end{enumerate}

\begin{thm}\label{homology of long knots in 3 space}
The homology of $\K_{3,1} $ satisfies the following properties.

\begin{enumerate}
    \item The rational homology of $\K_{3,1} $ is a free
    Gerstenhaber-Poisson
    algebra generated by $V = H_*(X_{\K}; \mathbb Q)$.
    \item The homology of $\K_{3,1}$ with $\mathbb F_p$ coefficients
    is a free restricted Gerstenhaber-Poisson algebra generated by
    $V = H_*(X_{\K}; \mathbb F_p)$ as described in \cite{CLM}.
    \item There are isomorphisms of Hopf algebras
\begin{enumerate}
    \item $S[\sigma^{-1}L[\sigma(V)]] \to H_*(\K_{3,1}; \mathbb Q)$
    for $V = H_*(X_{\K}; \mathbb Q)$,
    \item $S[\sigma^{-1}L^{(2)}[\sigma(V)]]\to H_*(\K_{3,1}; \F_2)$ for
    $V = H_*(X_{\K};\F_2)$, and
    \item $S[\sigma^{-1}L^{(p)}[\sigma(V)] \oplus \sigma^{-2}W^p[\sigma(V)]]
    \to H_*(\K_{3,1};\F_p)$
    for $V = H_*(X_{\K};\F_p)$ in case $p$ is an odd prime.
\end{enumerate}
These isomorphisms specialize to an identification of the homology
of each path-component of $\K_{3,1}$ with the one ambiguity that the
homology of the components of knots arising from hyperbolic
satellite operations is not given in a closed form here. More
information is described in Section \ref{htpe}.

\item The integer homology of $\K_{3,1}$ has $p$-torsion of
arbitrarily large order (with examples listed in Section
\ref{c2homol}, and \ref{khomol}).
\end{enumerate}

\end{thm}

A primary development in this paper is our recursive application of
the above theorem. 
Let $K_c \subset \K_{3,1}$ denote the subspace of $\K_{3,1}$ consisting of
all cable knots. There is a homotopy-equivalence $\Zed \times S^1
\times \K_{3,1} \to K_c$. Let $K_s \subset \K_{3,1}$ denote the
subspace of $\K_{3,1}$ which are connect-sums of any number of cable
knots.  Then there is a homotopy-equivalence $\Cu_2 (K_c \sqcup
\{*\}) \to K_s$. The composite of the two maps is a
homotopy-equivalence $\Cu_2 ( (\Zed \times S^1 \times \K_{3,1})
\sqcup \{*\} ) \to K_s$. Since $K_s$ is a collection of
path-components of $\K_{3,1}$, this map can be iterated, giving a
the homology of $\K_{3,1}$, as a Gerstenhaber-Poisson algebra, a
fractal-like structure. These statements will be justified in
Sections \ref{ops} and \ref{htpe}, and explored more fully in
Sections \ref{findtorsion} and \ref{tcns}.

These results lead to some natural questions about the structure of
the homology of the higher-dimensional embedding spaces $\K_{n,1}$
($n \geq 4$)  studied recently by Sinha \cite{Sinha}, Volic
\cite{Volic}, Lambrechts \cite{Lambrechts} as well as others
\cite{Altshuler.Freidel,CRL,Goodwillie.Weiss,K3,Sakai}.
Constructions related to these questions are also addressed here.

By Theorem \ref{homology of long knots in 3 space}, there is
arbitrarily large $p$-torsion in the homology of $\K_{3,1}$.
Examples of Theorem \ref{homology of long knots in 3 space} for
knots whose path-components have higher $2$-torsion in their integer
homology is given next. This higher torsion can be regarded as a
coarse ``measure" of the ``complexity" of a knot's
JSJ-decomposition.

\begin{enumerate}
    \item Let $\K_{3,1}(f)$ denote the path-component of a torus
knot $f$. Thus $\K_{3,1}(f)$ has the homotopy type of a circle
\cite{Hatcher4, B2}.
    \item Given any space $X$, and a strictly positive integer $q$, define
$$E(q,X) = \cfg(\Real^2,q) \times_{\Sigma_q}X^q.$$ Assume
that $\K_{3,1}(f)$ has the homotopy type of $S^1$.
$E(4,\K_{3,1}(f))$ has the property that

\begin{enumerate}
  \item $H_2( E(4,\K_{3,1}(f)) ) = \cyc{2}$,
  \item $H_3( E(4,\K_{3,1}(f))) = 0$,
  \item $H_3(S^1 \times E(4,\K_{3,1}(f)))$ is isomorphic to
  $\cyc{2}$.

\item Furthermore, $E(2^s, S^1 \times E(4,\K_{3,1}(f)))$ has the homotopy
    type of a path-component $\K_{3,1}(g)$ for a long knot $g$ as given in \cite{B},
    and has torsion of order $2^{s+1}$ in its integer homology by Section \ref{khomol},
    and \ref{findtorsion}.
    In this case, $g$ is a connected-sum of $2^s$ copies of the same summand, and
    that summand is a $p/q$-cable of a connected sum of four copies of the same torus knot.
    In the language of \cite{B1},
    $$g=\left( \left( T^{(p,q)} \rep_4 \splice \Keychain^4  \right)
    \splice \Seifert^{(p,q)} \right) \rep_{2^s} \splice \Keychain^{2^s}$$
    The elements of this notation is described in detail in
    \cite{B1} and is summarized in Section \ref{ops}.
    \item A second example is  $E(2^s, \K_{3,1}(h)))$ where
    $h=T^{(p,q)} \splice W$ and $W$ is the Whitehead link.
    In this case,
    $$\K_{3,1}(h) \simeq S^1 \times \left( S^1 \times_{\Sigma_2} S^1 \right)$$
    where $S^1 \times_{\Sigma_2} S^1$ is the Klein bottle \cite{B2}.
    $H_1( \K_{3,1}(h)) ) = \cyc{2} \oplus \Zed^2$.
\end{enumerate}
\end{enumerate}

A more complete description of the homology of $\K_{3,1}$ is given
in Sections \ref{c2homol}, and \ref{khomol}. The homology of each
path-component is given in terms of Theorem \ref{homology of long
knots in 3 space} as well as filtrations of the values of the
functors listed in that theorem.

Consider the subspace $\mathcal T\K_{3,1}$ of $\K_{3,1}$ consisting
of the union of the components of $\K_{3,1}$ corresponding to knots
which can be obtained from torus knots via iterations  of the
operations given by cabling and connected-sum. These are the knots
whose complements are `graph manifolds' ie: a union of
Seifert-fibered manifolds.  The structure of the homology of
$\mathcal T\K_{3,1}$ is described in Section \ref{tcns}. 
This is our primary source of $p^n$-torsion in $H_*(\K_{3,1};\Zed)$.

A further consequence of Theorem \ref{long knots in 3 space} is the
next result which follows directly, and is proven in Section
\ref{relations}. Let $\Emb_*(S^1,S^n)$ denote the space of smooth
pointed embeddings.

\begin{prop}\label{prop:embeddings of circles in the three sphere}
The group $SO(n-1)$ acts naturally on $\K_{n,1}$ (rotations that fix
the `long axis'), and there are morphisms of bundles for which
each vertical map is a homotopy equivalence:

\[
\begin{CD}
SO(n) \times_{SO(n-1)} \K_{n,1} @>{i}>> SO(n+1) \times_{SO(n-1)}
\K_{n,1} @>{p}>>
SO(n+1)/{SO(n)} \\
 @VV{\theta_n }V  @VV{\theta_n}V          @VV{1}V \\
Emb_*(S^1,S^n)  @>{j}>> \Emb(S^1,S^n)  @>{}>> S^n
\end{CD}
\]

Thus, there is a bundle  $$SO(n+1) \times_{SO(n-1)}\K_{n,1} \to
SO(n+1)/ SO(n-1)$$ with fibre $\K_{n,1}$. Furthermore, there is a
homeomorphism
$$\Emb(S^1,S^3) \to S^3 \times \Emb_*(S^1,S^3)$$ for which
$\Emb_*(S^1,S^3)$ denotes the space of smooth pointed embeddings,
and the bundle
$$\K_{3,1} \to \Emb_*(S^1,S^3) \to S^2$$
is the induced bundle with fibre $\K_{3,1}$ from the bundle
$$SO(2) \to SO(3) \to S^2.$$
\end{prop}

Section \ref{relations} gives precise relationships (such as the
above proposition) between the homotopy-type of the embedding spaces
$\K_{n,j}$, $\Emb(S^j,S^n)$ and $\Emb(S^j,\Real^n)$.

Notice that the homological properties of each path component thus
give knot invariants. This is illustrated by the following
proposition.

\begin{prop}\label{prop:Two torsion in h1}
\begin{enumerate}
    \item A knot $f : S^1 \to S^3$ is the unknot if and only if
its component in $\Emb(S^1,S^3)$ contains no $2$-torsion in its 1st
homology group.
    \item A long knot $f : \Real^1 \to \Real^3$ in $\K_{3,1}$
    is the `long unknot' if and only if its component has trivial
    first homology group.
    \item An embedding of $S^1$ in $\Real^3$ is the unknot if and only
if its component in $\Emb(S^1,\Real^3)$ has torsion first homology group.  
It is also true if and only if its 2nd homology group is trivial.
\end{enumerate}
\end{prop}

\begin{thm}\label{H1.of.the knot.space}
Let $\K_{3,1}(f)$ denote a path-component of $\K_{3,1}$. Then
$$H_1 (\K_{3,1}(f);\Zed)$$ is a finite direct-sum of copies of
$\Zed$ and $\cyc{2}$.
\end{thm}

In addition, a characterization of the components of $\K_{3,1}$ such
that $H_1(\K_{3,1}(f);\Zed)$ contains $2$-torsion is given in
Section \ref{initial_computations_of_H}. A precise identification of
those knots $f$ such that $H_1(\K_{3,1}(f);\Zed)$ contains a
$\cyc{2}$ summand is also given. In Section \ref{no_p_torsion} the
least degree in which odd $p$ torsion in $H_*(\K_{3,1};\Zed)$ occurs
is as follow.

\begin{thm}\label{first.odd.p.torsion}
Let $f$ denote a long knot and $p$ an odd prime. If
$H_i(\K_{3,1}(f);\Zed)$ contains $\cyc{p}$, then $i \geq 2p-2$.
\end{thm}

Much recent progress has been made on the structure of spaces of
embeddings via finite-dimensional model spaces and approximations.
Some of this was first given by Vassiliev \cite{Vass} and has been
the subject of further study via the Goodwillie Calculus of
Embeddings by Sinha \cite{Sinha}, Volic \cite{Volic} Lambrechts
\cite{Lambrechts}, or cohomological techniques such as Bott-Taubes
integrals \cite{Bott-Taubes, CRL}.

The Gerstenhaber-Poisson algebra above was first considered on the
$E^2$-level of the Vassiliev spectral sequence by Tourtchine
\cite{Tourtchine}. Related progress is given in work of
Altshuler-Freidel \cite{Altshuler.Freidel}, Bar-Natan \cite{BN},
Cattaneo, Cotta-Ramusino-Longoni \cite{CRL}, Kohno \cite{K3},
Kontsevich \cite{Kont}, Lescop \cite{Lescop}, Polyak-Viro
\cite{PolyViro}, Sakai \cite{Sakai}, Watanabe \cite{Watan} as well
as others.

This paper takes the direction of using Gramain and Hatcher's
techniques for understanding the homotopy type of $\K_{3,1}$, one
component at a time \cite{GramainPi1, Hatcher4}. The central
construction of Hatcher \cite{Hatcher4} is to consider the
components of the knot space as the classifying space of the mapping
class group of the knot complement. One then studies how such a
mapping class group acts on the JSJ-tree of the knot complement as
in \cite{B, B2}, using Hatcher's results on the homotopy type of
diffeomorphism groups of Haken manifolds \cite{Hatcher1} to assemble
an answer. Thus most of the results here are complementary to the
results of the authors mentioned in the previous two paragraphs.

The authors would like to thank the University of Tokyo, the Max
Planck Institute for Mathematics in Bonn, Institut des Hautes 
\'Etudes Scientifiques, the Institute for
Advanced Study, the Pacific Institute of Mathematics and
the American Institute of Mathematics for partial
support during the preparation of this paper.

{\bf TABLE OF CONTENTS}

\begin{description}
    \item[1] Introduction
    \item[2] Notation, labelling components
    \item[3] The homotopy type of $\K_{3,1}$
    \item[4] Relations among various spaces
    \item[5] On the homology of $\Cu_2(X \amalg \{*\})$
    \item[6] On the homology of $\K_{3,1}$
    \item[7] Higher $p$-torsion in the integer homology of $\K_{3,1}$
    \item[8] $H_1 (\K_{3,1};\Zed)$
    \item[9] The first appearance of odd $p$-torsion
    \item[10] On the subspace generated by torus knots
    \item[11]Closed knots and homology
    \item[12] Problems
\end{description}

\section{Notation, labelling components}\label{ops}

Whitney \cite{Whitney3} showed that the embedding space $\K_{n,j}$
is connected for $n > 2j+1$.  By work of Wu \cite{Wu}, $\K_{n,j}$ is
also connected provided both $n > 2j \text{ and } j>1$. That
$\K_{n,1}$ is also connected for $n=1$ is elementary. The fact that
$\K_{2,1}$ is connected is equivalent to the smooth
Alexander/Schoenflies theorem in dimension 2. In co-dimension $3$
and higher. Haefliger \cite{Haf2} vastly generalized Whitney's
result, proving that $\K_{n,j}$ is connected provided $2n > 3j+3$,
and $\pi_0 \K_{n,j}$ is non-trivial for $2n=3j+3$. This work has
recently been extended by the first author to a computation of the
first non-trivial homotopy group of $\K_{n,j}$ provided $2n-3j-3
\geq 0$ \cite{B3}.

When $2n \leq 3j + 3$ the space $\K_{n,j}$ could potentially have
many connected components. $\pi_0 \K_{n,j}$ was shown to be a
group by Haefliger \cite{Haf2} provided $n-j>2$, whereas it is
only a monoid for $n-j \leq 2$.  A fundamental example is the space
$\K_{3,1}$ which has countably infinite many components, and no
inverses in the monoid $\pi_0 \K_{3,1}$ \cite{Sch2}. Given $f \in \K_{n,j}$,
let $\K_{n,j}(f)$ denote the path-component of $\K_{n,j}$ containing
$f$.

We will use the notation $\EK{1,D^{n-1}}$ as defined in \cite{B} for
the space of framed long knots in $\Real^n$. Given a compact
manifold $M$, define
$$\EK{k,M} = \{ f \in \EC{k}{M}, supp(f) \subset \I^k \times M\}.$$
 Here the support of
$f$, $supp(f)$ is defined by $supp(f) = \{ x \in \Real^k \times M :
f(x)\neq x \}$ and $\I = [-1,1]$. 
$\EK{1,D^{n-1}}$ is not homotopy equivalent to $\K_{n,1}$ in
general, but as described in \cite{B} there is a fibration
$$\Omega SO(n-1) \to \EK{1,D^{n-1}} \to \K_{n,1}$$
which splits at the fibre (via a $2$-cubes map) for $n\in
\{1,2,3\}$, allowing us to think of $\K_{3,1}$ as a sub $2$-cubes object of
$\EK{1,D^2}$.

This section collects information on the indexing of the components
of $\K_{3,1}$ which is given in terms of the `companionship tree'
classification of knots, an application of the Jaco-Shalen-Johannson
(JSJ) decomposition of knot complements. The indexing that we will
use is described in detail in \cite{B1}. Aspects of this indexing
have been partially described before in the works of Budney
\cite{B}, Eisenbud and Neumann \cite{EN}, Schubert \cite{Sch2}, and
the unpublished work of Bonahon and Siebenmann \cite{BS}, as well as
the survey work of Kawauchi \cite{Kaw}. Indeed, the results in
\cite{B1} should be thought of as a uniqueness statement for
Schubert's satellite operations that he describes in \cite{Sch2}. In
the book by Eisenbud-Neumann \cite{EN} this method of indexing is
called the splice decomposition of links, but is specialized to the
case of links in homology spheres whose complements are graph
manifolds. A terse statement of the results in \cite{B1} given next
suffice for the applications here. More complete as well as more
specific information is given in \cite{B1}.

\begin{defn}
An $n$-component link in $S^3$ is a compact, connected, oriented,
$1$-dimensional submanifold of $S^3$ consisting of $n$ path
components labelled with distinct numbers from the set
$\{0,1,2,\cdots,n-1\}$. Thus the notation
$L=(L_0,L_1,\cdots,L_{n-1})$ is used frequently for $n$-component
links. A knot $K$ (in $S^3$) is a $1$-component link.

An $n$-component link $L$ is the unlink if there exists $n$
disjointly embedded $2$-discs $D=(D_0,D_1,\cdots,D_{n-1})$ in $S^3$
whose boundary is $L$, $\partial D = (\partial D_0,\partial D_1,
\cdots , \partial D_{n-1})=(L_0,L_1,\cdots,L_{n-1})=L$.

For $n \geq 0$ an $(n+1)$-component link $L=(L_0,L_1,\cdots,L_n)$ is
said to be a KGL (knot generating link) if the sublink
$(L_1,L_2,\cdots,L_n)$ is the unlink.

Given an $(n+1)$-component KGL $L$ and $n$ knots
$J=(J_1,J_2,\cdots,J_n)$ in $S^3$ there is an operation called
splicing defined in \cite{B1} which produces a knot $J \splice L$ in
$S^3$.  Here is a rough statement of the splicing construction. Fix
$D=(D_1,\cdots,D_n)$ $n$ disjointly embedded discs in $S^3$ such
that $\partial D = (L_1,L_2,\cdots,L_n)$. Let 
$\nu_D : [-\infty,\infty] \times D^2 \to S^3$ be a closed
tubular neighbourhood of $D$. Let $C_{L'}$ be the complement of an open
tubular neighbourhood of $(L_1,\cdots,L_n)$ in $S^3$, and define
$R : C_{L'} \to S^3$ to be unique continuous function which is the
identity outside $img(\nu D)$, and on the image of $\nu_D$ define it
to be the conjugate $\nu_D \circ (\tilde J_1 \sqcup \cdots \tilde J_n) \circ \nu_D^{-1}$,
where $\tilde J_i \in \EK{1,D^2}$ is the framed long knot in the homotopy-fibre
of the map $\EK{1,D^2} \to \Omega SO(2)$ corresponding to $J_i$ under the map
$\EK{1,D^2} \to \K_{3,1} \to \Emb(S^1,S^3)$. $J \splice L$ is defined to be
the image of $L_0$ under the embedding $R$. See \cite{B1} for details.
\end{defn}

\begin{eg}
Let $W$ denote the Whitehead link and $F_8$ the figure-8 knot. 
\begin{figure}[htp]
\centering
{
\psfrag{W}[tl][tl][1][0]{$W$}
\psfrag{F}[tl][tl][1][0]{$F_8$} 
\psfrag{n}[tl][tl][1][0]{$K=F_8\splice W$}
$$\includegraphics[width=7cm]{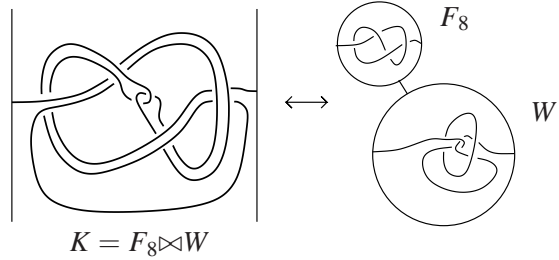}$$
}
\caption{Whitehead double and companionship tree}\label{fig:f1}
\end{figure}
\end{eg}

The role of splicing is that it is an operation that takes knots and
$KGL$'s as input and produces a knot of greater `complexity' in the
sense that the companionship trees of the input data is spliced
together to produce the companionship tree of $J \splice L$.  We
proceed to make these ideas more precise.

\begin{defn}
The Hopf link $\Keychain^1$ is the $2$-component link in $S^3$ given
by
$$\{ (z_1,0) \in \Complex^2 : z_1 \in \Complex, |z_1|=1 \} \cup \{ (0,z_2) : z_2 \in \Complex, |z_2|=1 \} \subset S^3$$
where $S^3$ is regarded the unit sphere in  $\Complex^2$. 

If one takes a connected-sum of $p$ copies of the Hopf link along a common
component, one obtains the $(p+1)$-component link, which we will call the
$(p+1)$-component keychain link $\Keychain^p$ (see Figure \ref{fig:f2}).
$$\{ (z_1,0) \in \Complex^2 : |z_1|=1\} \cup \bigcup_{k=1}^p
  \{ \frac{1}{\sqrt{2}}(e^{ \frac{2\pi ik}{p}}, z_2) : |z_2|=1 \} \subset S^3$$

\begin{figure}[htp]
\centering
{ 
\psfrag{...}[tl][tl][0.8][0]{$\cdots$}
\psfrag{h1}[tl][tl][1][0]{$\Keychain^1$}
\psfrag{hopf}[tl][tl][1][0]{$\Keychain^p$}
\includegraphics[width=2cm]{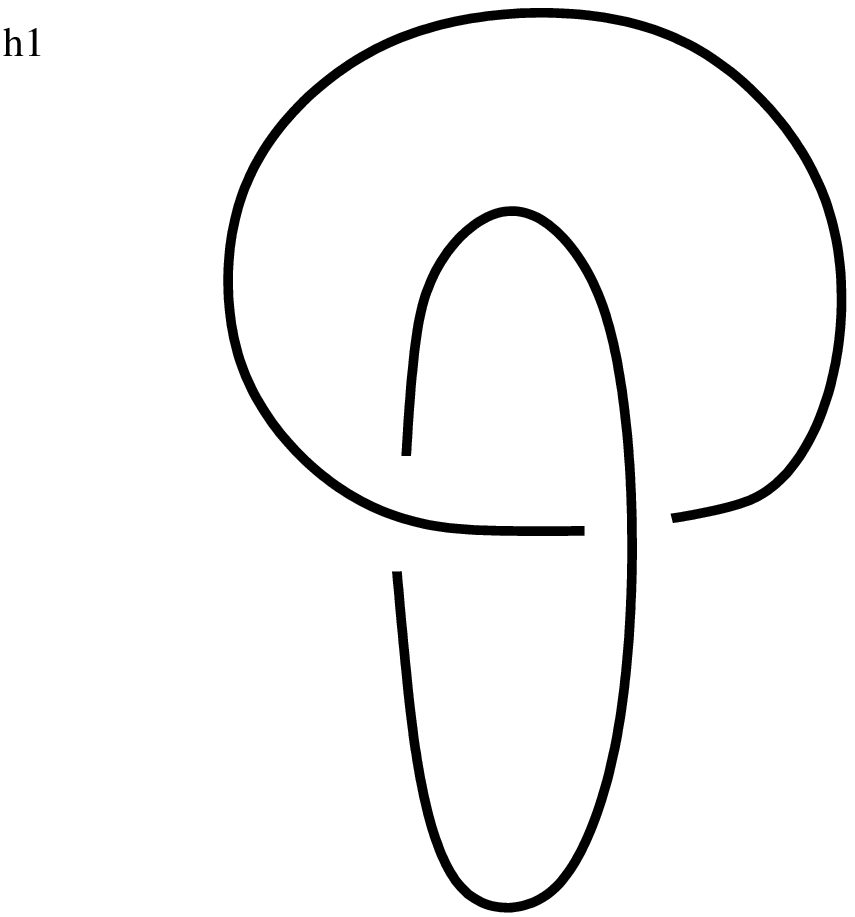}\hskip 10mm \includegraphics[width=4cm]{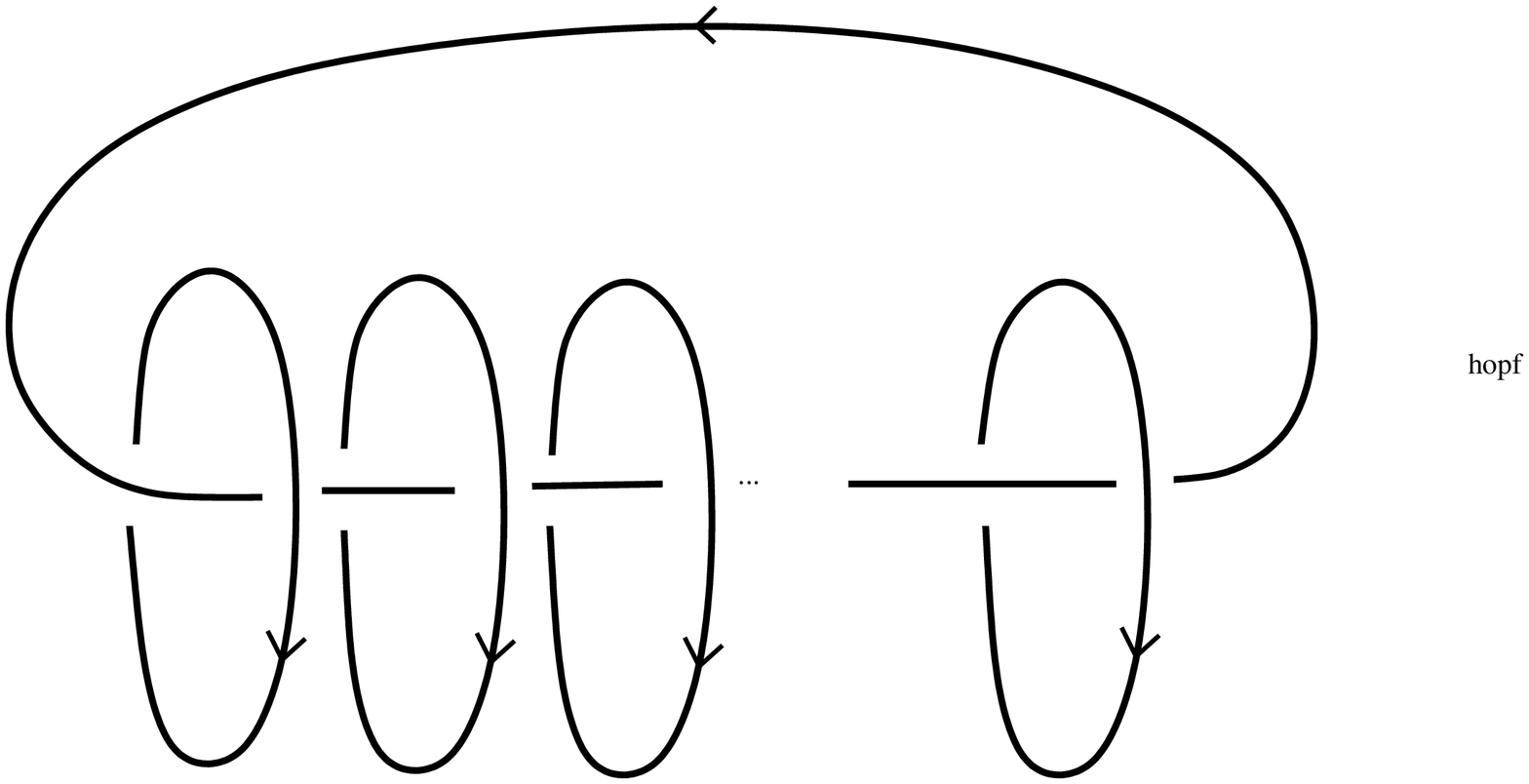}
}
\caption{Hopf link and keychain link}\label{fig:f2}
\end{figure}

For any $(p,q) \in \Zed\times \Natural$, the $(p,q)$-Seifert link
$\Seifert^{(p,q)}$ is defined to be
$$\{ (z_1,0) \in \Complex^2 : |z_1|=1 \} \cup \{ (z_1,z_2) \in \Complex^2 : |z_1|=|z_2|=\frac{1}{\sqrt{2}}, z_1^q=z_2^p \} \subset S^3 $$
The $(p,q)$-Seifert link has $GCD(p,q)+1$ components (see Figure \ref{fig:f3}).

\begin{figure}[htp]
\centering
{ 
\psfrag{p}[tl][tl][1][0]{$p$} \psfrag{q}[tl][tl][1][0]{$q$}
\psfrag{b}[tl][tl][1][0]{$\Seifert^{(p,q)}$}
\includegraphics[width=5cm]{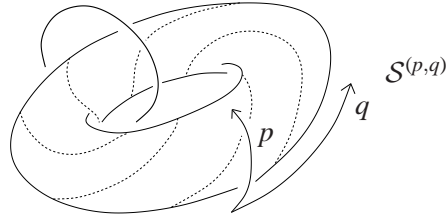}
}
\caption{Seifert link}\label{fig:f3}
\end{figure}

For any $(p,q) \in \Zed \times \Natural$, $GCD(p,q)=1$, the
$(p,q)$-torus knot $T^{(p,q)}$ is
$$\{ (z_1,z_2) \in \Complex^2 : |z_1|=|z_2|=\frac{1}{\sqrt{2}}, z_1^q=z_2^p \} \subset S^3 $$
\end{defn}

\begin{thm} \cite{B1}
Given a knot $K$ in $S^3$ there is a finite, labelled, rooted
tree-valued invariant of $K$ denoted $\IG_K$ having the following
properties:
\begin{enumerate}
\item Each vertex of the tree is labelled by a link and any link from the following list
is admissible:
 \begin{enumerate}
 \item Torus knots $T^{(p,q)}$ for $p/q \in \Rational$, $GCD(p,q)=1$, $q \geq 2$.
 \item Seifert links $\Seifert^{(p,q)}$ for $GCD(p,q)=1$, $q \geq 1$.
 \item Keychain links $\Keychain^p$ for $p \geq 2$.
 \item Hyperbolic KGLs.
 \item The unknot.
 \end{enumerate}
\item Given any vertex in $\IG_K$, the number of children of the vertex is one less than the
number of components of the link that decorates the vertex.
\item If any vertex is decorated by a keychain link $\Keychain^p$, none of its children are allowed to
be decorated by keychain links.
\item A vertex of the tree $\IG_K$ can be decorated by the unknot if and only if the tree $\IG_K$
consists of only one vertex.
\item If one changes all the labels on the tree $\IG_K$
by substituting for each vertex label $L$ its complement $C_L$ one
obtains $G_K$, the JSJ-tree of $K$ \cite{B1}. This is the tree whose
vertex set is the set of path components of the knot complement
$C_K$ split along its JSJ-tori, and the edges are the JSJ-tori of
$C_K$.
\item If $\IG_K$ consists of more than one vertex, then $K = J \splice L$ where
the root of $\IG_K$ is labelled by $L$ and $\IG_{J_i}$ are the
subtrees rooted at the children of $L$ in $\IG_K$.
\item The number of vertices of $\IG_K$ is one more than the number of
tori in the JSJ-decomposition of the complement of $K$ in $S^3$.
Thus for example, $\IG_K$ is a one-vertex tree if and only if $K$ is
either hyperbolic, a torus knot, or the unknot.
\end{enumerate}

The above properties 1 through 4 are complete, in the sense that any
tree satisfying properties 1 through 4 is realizable as $\IG_K$ for
some knot $K$. $\IG_K$ is known as the companionship tree of $K$.

Given a vertex $v$ of $\IG_K$, there is a maximal subtree of $\IG_K$
rooted at $v$, and this subtree is the companionship tree of a
unique knot in $S^3$, $K_v$.  $K_v$ is called a companion knot to
$K$.

Item (6) implies that if one writes down the `postorder' (reverse
Polish) listing of $\IG_K$, one is simply writing $K$ as an iterated
splice knot where all the KGL's used in splicing come from the list
(1). Thus $\IG_K$ could simply be considered a precise way to
specify $K$ as a splice of atoroidal KGLs.
\end{thm}

Unlike links with Seifert-fibred complements, hyperbolic KGLs have
no known canonical enumeration.

Some elementary examples of hyperbolic KGLs are:
\begin{itemize}
\item the figure-8 knot
\item the Whitehead link
\item the Borromean rings
\end{itemize}
Hyperbolic KGLs of arbitrarily many components are known to exist by
the work of Kanenobu \cite{Kan}. For details on the hyperbolic
structures, see for example the textbook of Thurston \cite{Thur}.

\begin{figure}[htp]
\centering
\centerline{ { \psfrag{torusknot}[tl][tl][1][0]{A trefoil knot $K=T^{(-3,2)}$} 
               \psfrag{blah}[tl][tl][1][0]{\ \ \ \ \ \ \ \ \ $\IG_K = K$} 
               \includegraphics[width=3cm]{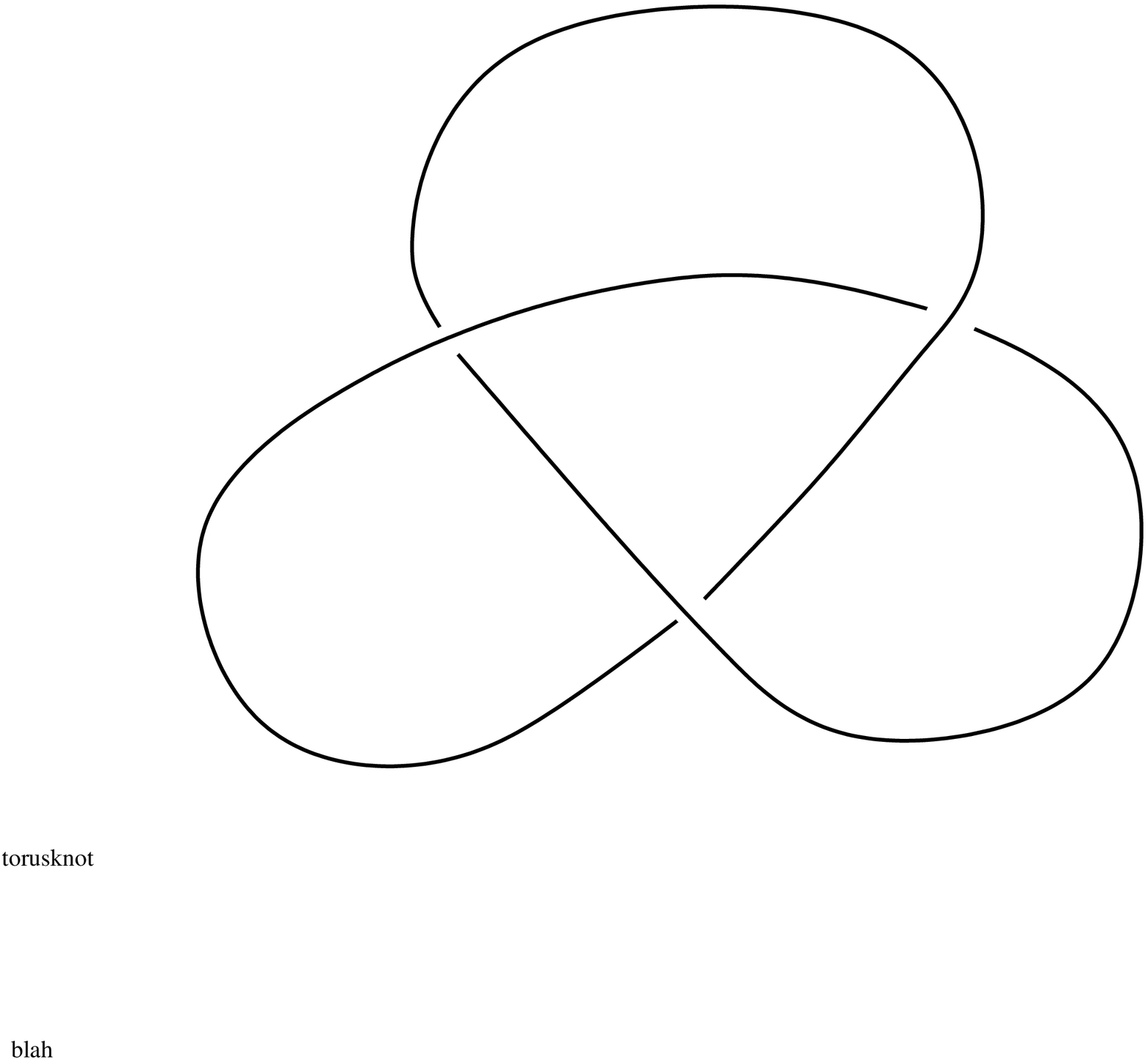} \hskip 3cm
               \psfrag{cableknot}[tl][tl][1][0]{ $T^{(-3,2)} \splice \Seifert^{(17,2)}$} 
               \includegraphics[width=3cm]{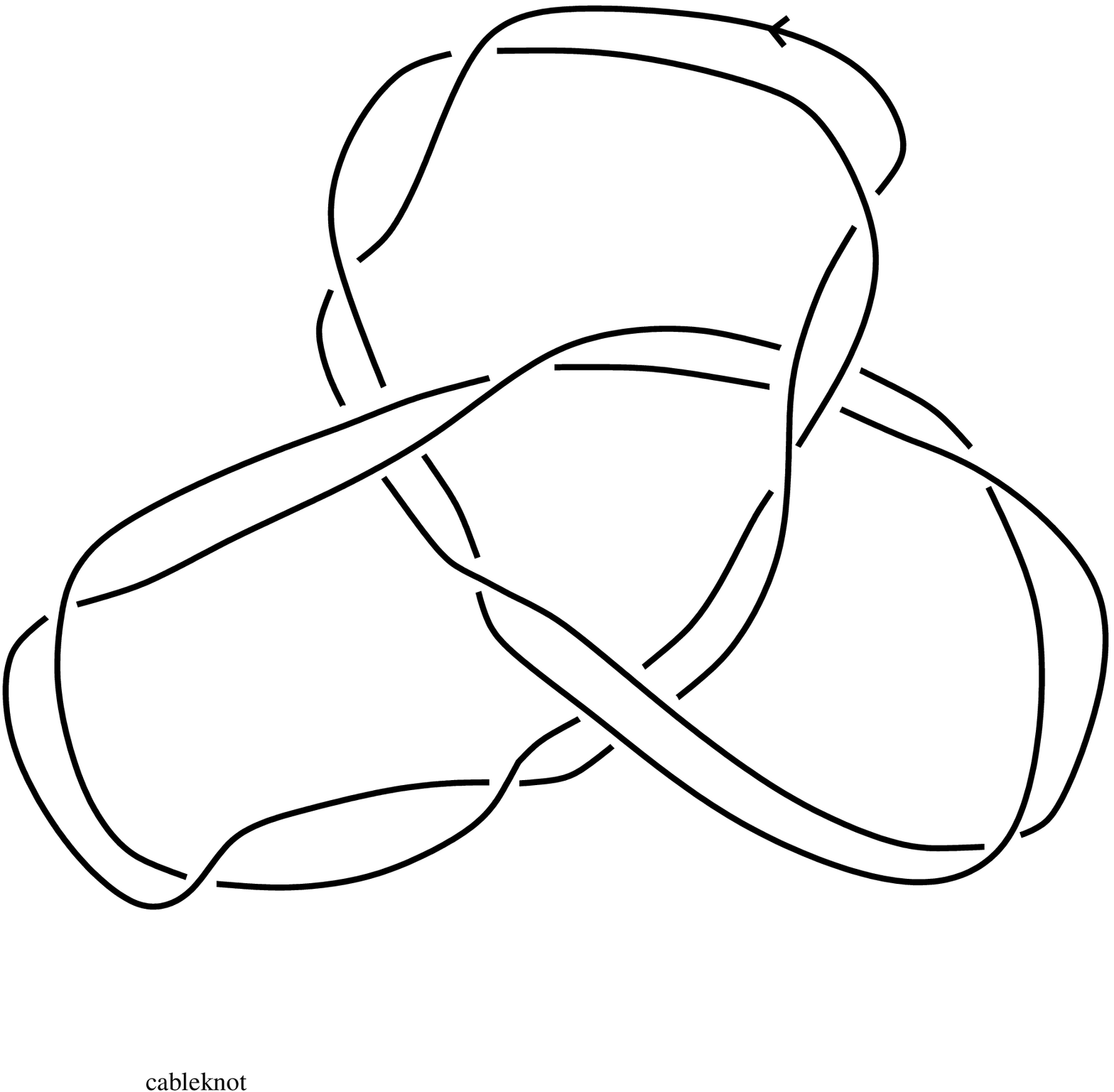} } }
\caption{Trefoil and (17,2)-cable}\label{fig:f4}
\end{figure}

\begin{figure}[htp]
\centering
\centerline{ { \psfrag{n}[tl][tl][1][0]{$F_8$}
               \includegraphics[width=3.3cm]{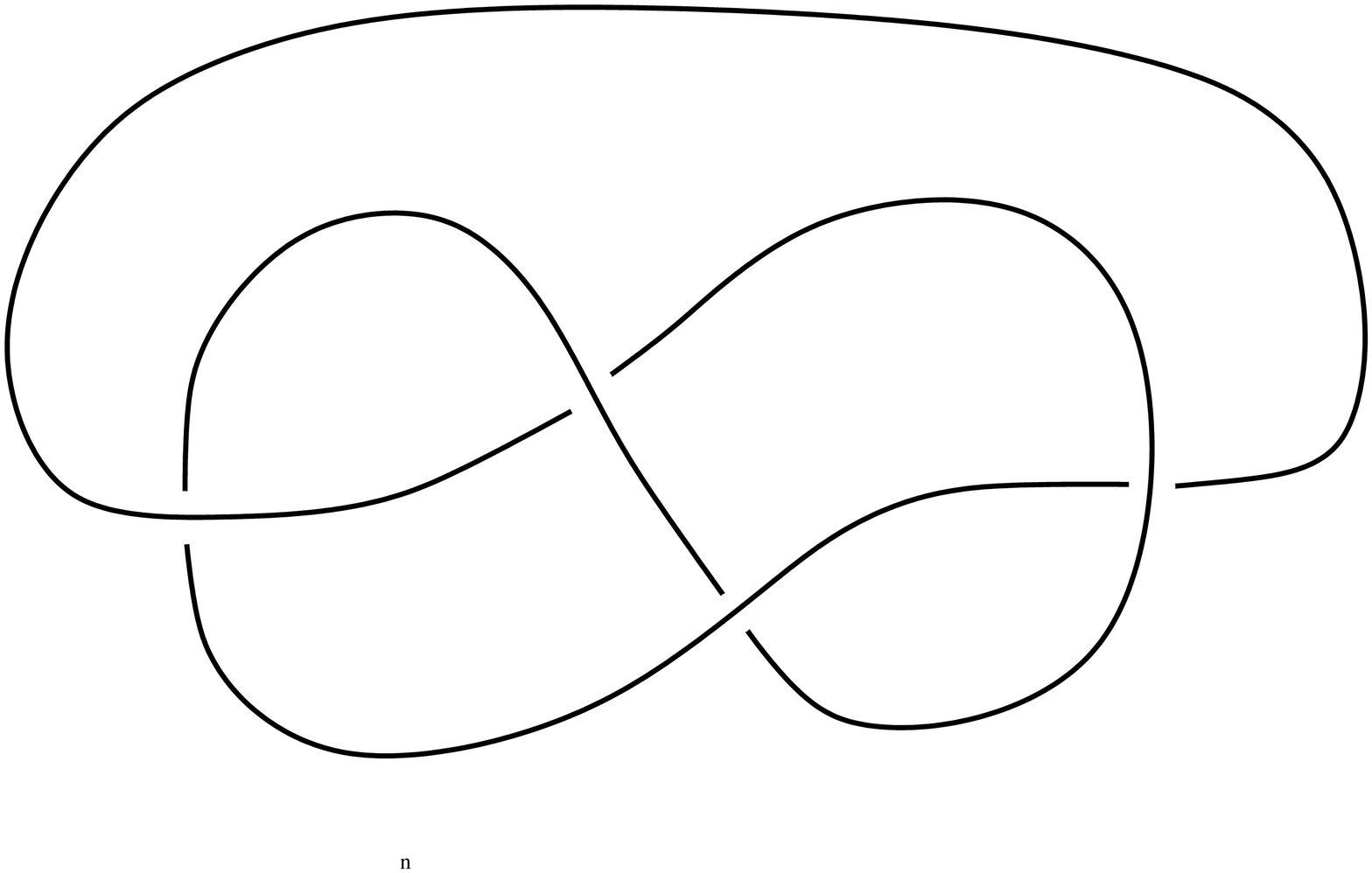} \hskip 1cm 
               \psfrag{cs}[tl][tl][1][0]{$((T^{(3,2)}, F_8) \splice \Keychain^2) \splice \Seifert^{(-17,2)} $}
               \includegraphics[width=6cm]{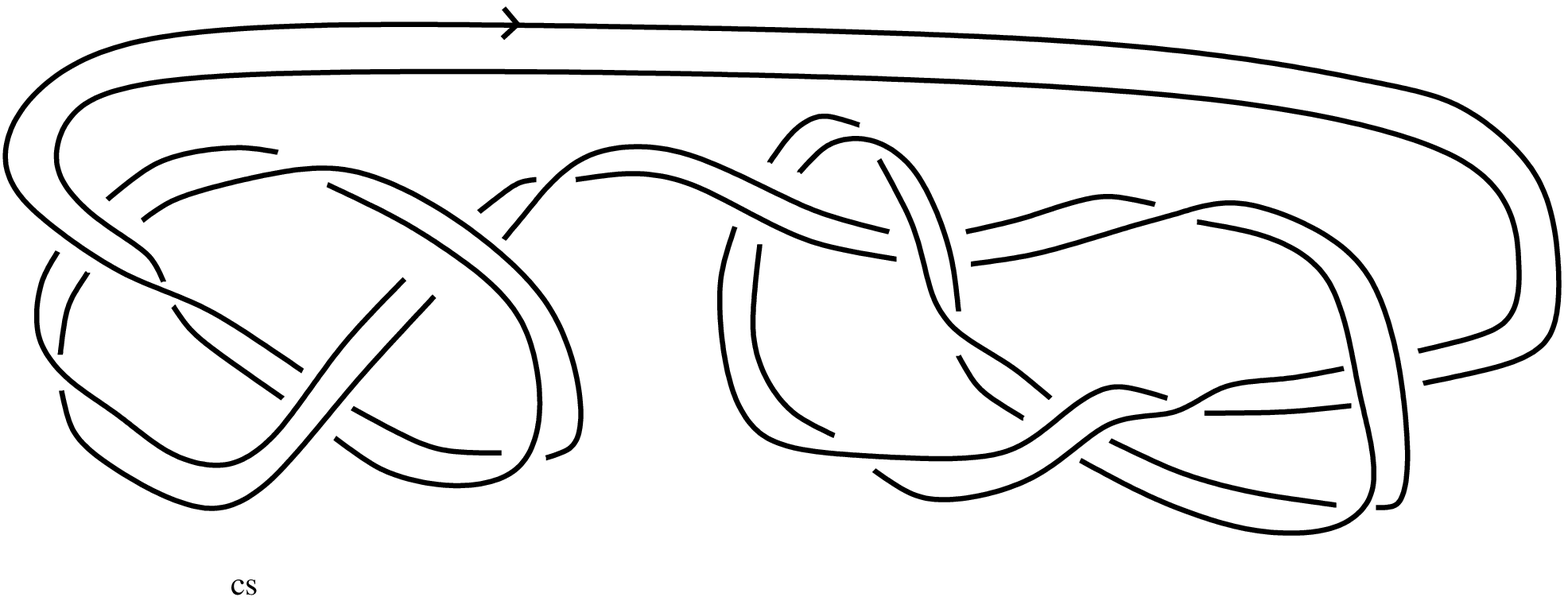} } } 
\caption{Cable of connect sum of trefoil and figure-8 knot}\label{fig:f5}
\end{figure}

Figures \ref{fig:f4} and \ref{fig:f5} give examples of knots $K$ and 
the associated tree $\IG_K$, and the corresponding splice notation, 
where $F_8$ denotes the figure-8 knot.
Let $B=(B_0,B_1,B_2)$ denote the Borromean rings, and let $B_{i,j}$
be the $3$-component link in $S^3$ obtained from $B$ by doing $i$
Dehn twists about the spanning disc of $B_1$ and $j$ Dehn twists
about the spanning disc for $B_2$.

\begin{figure}[htp]
\centering
\centerline{ { \psfrag{K}[tl][tl][1][0]{$(F_8,T^{(3,2)}) \splice B$}
               \psfrag{IG}[tl][tl][1][0]{$\IG_K$}
               \includegraphics[width=6cm]{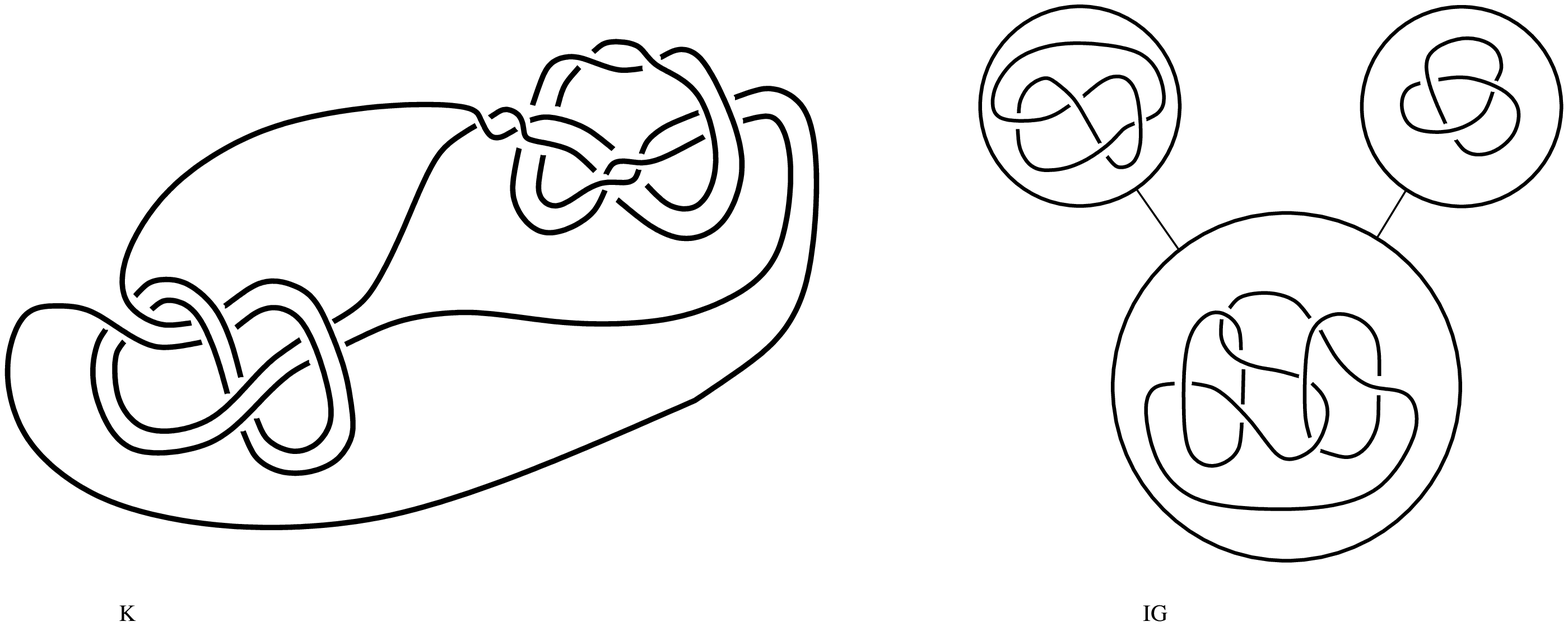} \hskip 0.5cm
               \psfrag{S2}[tl][tl][1][0]{$\IG_K$}
               \psfrag{K2}[tl][tl][1][0]{$(F_8,T^{(3,2)}) \splice B_{0,3}$}
               \includegraphics[width=6cm]{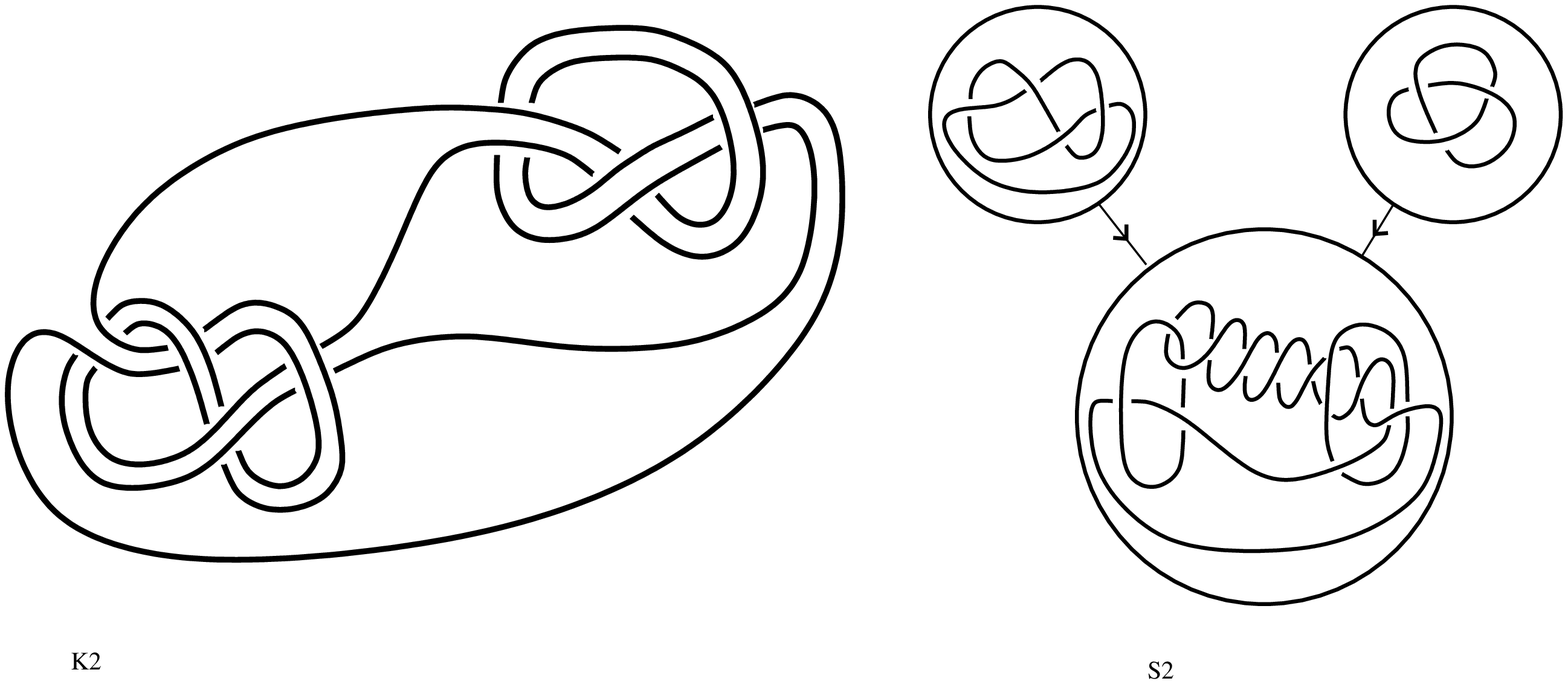} } }
\caption{Various Borromean splices} \label{fig:f6}
\end{figure}

The spaces $\EK{1,D^{n-1}}$ admit
an action of the operad of little $2$-cubes \cite{B}.
 Using the connectedness of
$\K_{n,1}$ for $n \geq 4$ together with the cubes action one can
prove that $\EK{1,D^{n-1}}$ has the homotopy type of a $2$-fold loop
space for $n \geq 4$. At present it is not known what the 2-fold
de-looping of $\EK{1,D^{n-1}}$ is.  Recently, P.~Salvatore
\cite{Sal2} has constructed an action of the operad of $2$-cubes on
$\K_{n,1}$ for all $n \geq 4$.

The previously described fibre bundle $$ \Emb(S^1,S^3) \to
SO(4)/SO(2)$$ whose fibre inclusion $i:\K_{3,1} \to \Emb(S^1,S^3)$
is induced by the one-point compactification
is explored more deeply in Section \ref{relations}. For the purpose
of this section and the study of the components of $\K_{3,1}$ and
$\Emb(S^1,S^3)$ respectively, we note that the inclusion $\K_{3,1}
\to \Emb(S^1,S^3)$ induces bijection on path-components. Thus, our
indexing of $\pi_0 \Emb(S^1,S^3)$ above by companionship trees
$\IG_K$ is also an indexing of $\pi_0 \K_{3,1}$.

\section{The homotopy type of $\K_{3,1}$}\label{htpe}

A detailed description of the homotopy type of $\K_{3,1}$ is given
in this section. This description is given in terms of the splicing
operations as described in Section \ref{ops}. A good general
reference for these results is the work \cite{B2}.

\begin{enumerate}
\item If $f$ is the unknot then $\K_{3,1}(f)$ is contractible by
work of Hatcher \cite{Hatcher4}.

\item If $f$ is a $p/q$-cabling of $g$ then by work of Hatcher
\cite{Hatcher4}, there is a homotopy equivalence
$$S^1 \times \K_{3,1}(g) \to \K_{3,1}(f).$$
We consider a torus knot to be a cable of the unknot, so we are
claiming all non-trivial torus knots $f$ satisfy $\K_{3,1}(f) \simeq
S^1$.

\item If $f = (f_1,f_2,\cdots,f_n) \splice \Keychain^n$
where $\{f_i : i\in \{1,2,\cdots,n\}\}$ are the prime summands of
$f$ and $n \geq 2$, then there is a homotopy equivalence
$$\Cu_2(n) \times_{\Sigma_f} \prod_{i=1}^n \K_{3,1}({f_i}) \to \K_{3,1}(f)$$
where $\Sigma_f \subset \Sigma_n$ is the Young subgroup corresponding to
the partition $\sim$ of $\{1,2,\cdots,n\}$ given by $i \sim j
\Leftrightarrow \K_{3,1}(f_i)=\K_{3,1}(f_j)$. This result originally
appears in \cite{B}.

\item If a knot $f = (f_1,f_2,\cdots,f_n)\splice L$
where $L$ is a hyperbolic KGL then there is a homotopy-equivalence:

$$S^1 \times \left( SO(2) \times_{A_f} \prod_{i=1}^n \K_{3,1}(f_i) \right) \to \K_{3,1}(f).$$

\noindent we define $A_f$ as a subgroup of $B_L$. $B_L$ is the
subgroup of the group of hyperbolic isometries of the complement of
$L$ in $S^3$ which:
\begin{itemize}
\item extend to diffeomorphisms of $S^3$.
\item the extensions preserve $L_0$ and its orientation (ie: they
act on the cusp corresponding to $L_0$ by translations).
\item put together, the above two properties imply there is
an embedding $B_L \to \Diff(S^3;L,L_0)$.
\end{itemize}
It is a non-trivial fact \cite{B2} that the composite is an
embedding of groups $B_L \to \Diff(S^3;L,L_0) \to \Diff^+(L_0)$
where $\Diff^+(L_0)$ is the group of orientation-preserving
diffeomorphisms of $L_0$.   
Regard $B_L$ as a finite subgroup of $SO(2)$. There is a
representation of $B_L$ given by the composite: 
$B_L \to \Diff(S^3,L,L_0) \to \Diff(\sqcup_{i=1}^n L_i) \to \pi_0
\Diff(\sqcup_{i=1}^n L_i) \equiv \Sigma_n^+$ where we identify $\pi_0
\Diff(\sqcup_{i=1}^n L_i)$ with $\Sigma^+_n$, the signed symmetric group on
$\{1,2,\cdots,n\}$. $\Sigma_n$ acts on $\K_{3,1}^n$ by permutation of factors. 
$\Sigma_2$ acts on $\K_{3,1}$ by knot inversion -- fix an axis perpendicular
to the long axis, and rotate a knot by $\pi$ about this axis, this
is knot inversion. Stated another way, the group of rotations which
preserve the long axis, $O(2) \subset SO(3)$, acts on $\K_{3,1}$ by
conjugation. Fix an element $\varpi \in O(2) \setminus SO(2)$, then
$\varpi$ acts as an involution on $\K_{3,1}$, thus defining an
action of $\Sigma_2$ on $\K_{3,1}$. These two actions extend to an
action of $\Sigma_n^+$ on $\K_{3,1}^n$. $A_f$ is the subgroup of $B_L$
that preserves the path-component $\prod_{i=1}^n \K_{3,1}(f_i)$ of
$\K_{3,1}^n$.
\end{enumerate}

As mentioned in part (4) above, $\K_{3,1}$ is naturally an
 $O(2)$-space. Parts (1), (2) and (3) above
all are $O(2)$-equivariant homotopy equivalences, as shown in
\cite{B2}. Case (4) is only an $SO(2)$-equivariant
homotopy-equivalence, although the homotopy-class of $\varpi$ acting
on $\K_{3,1}$ is computed.

\section{Relations among various spaces}\label{relations}

The goal of this section is to compare the homotopy types of the
spaces:

\begin{itemize}
\item $\K_{n,1}$
\item $\Emb(S^1,S^n)$
\item $\Emb(S^1,\Real^n)$
\end{itemize}

The space of pointed, smooth embeddings $\Emb_*(S^1,S^n)$ will be a
useful auxiliary space.  Relationships between the embedding spaces
$\K_{n,j}$, $\Emb(S^j,S^n)$ and $\Emb(S^j, \Real^n)$ will also be
listed.

\begin{prop}\label{closed-long-prop} For all $n \geq 1$ there are morphisms of
fibrations for which each vertical map is a homotopy equivalence:

\[
\begin{CD}
SO(n) \times_{SO(n-1)} \K_{n,1} @>{i}>> SO(n+1) \times_{SO(n-1)}
\K_{n,1} @>{p}>>
SO(n+1)/{SO(n)} \\
 @VV{\Theta_n }V  @VV{\Theta_n}V          @VV{1}V \\
Emb_*(S^1,S^n)  @>{j}>> \Emb(S^1,S^n)  @>{}>> S^n
\end{CD}
\]

\[
\begin{CD}
\K_{n,1} @>{i}>> SO(n) \times_{SO(n-1)} \K_{n,1} @>{p}>>
SO(n)/{SO(n-1)} \\
 @VV{1}V  @VV{\Theta_n}V          @VV{1}V \\
\K_{n,1}  @>{j}>> \Emb_*(S^1,S^n)  @>{}>> S^{n-1}.
\end{CD}
\]

\end{prop}

\begin{proof}
Consider the maps $\Theta_n:SO(n+1) \times \K_{n,1} \to\
\Emb(S^1,S^n)$ obtained from the natural $SO(n+1)$-action on
$\Emb(S^1,S^n)$ together with the natural inclusion $\K_{n,1} \to\
\Emb(S^1,S^n)$. Notice that the map $\Theta_n$ is
$SO(n-1)$-equivariant and thus there is an induced map
$$\Theta_n: SO(n+1) \times_{SO(n-1)} \K_{n,1} \to\ \Emb(S^1,S^n).$$

Consider the natural fibrations
$$\Emb_*(S^1,S^n) \to\ \Emb(S^1,S^n) \to S^n, $$ and
$$\K_{n,1} \to\ \Emb_*(S^1,S^n) \to S^{n-1}.$$
The first map is a fibration by the isotopy extension theorem.
Indeed, Palais proved that in general `restriction maps' are locally
trivial fibre bundles \cite{Pal}. The map $\Emb_*(S^1,S^n) \to
S^{n-1}$ is the composition of the restriction map $\Emb_*(S^1,S^n)
\to \Emb_*(U,S^n)$ with the homotopy equivalence $\Emb_*(U,S^n) \to
S^{n-1}$ given by the derivative at $*$ where $U$ is some closed
interval neighbourhood of $*$ in $S^1$. Thus there is a map of
fibrations:

\[
\begin{CD}
SO(n) \times_{SO(n-1)} \K_{n,1} @>{i}>> SO(n+1) \times_{SO(n-1)}
\K_{n,1} @>{p}>> SO(n+1)/{SO(n)} \\
 @VV{\Theta_n }V  @VV{\Theta_n}V          @VV{1}V \\
Emb_*(S^1,S^n)  @>{j}>> \Emb(S^1,S^n)  @>{}>> S^n
\end{CD}
\] as well as

\[
\begin{CD}
\K_{n,1} @>{i}>> SO(n) \times_{SO(n-1)} \K_{n,1} @>{p}>>
SO(n)/{SO(n-1)} \\
 @VV{1}V  @VV{\Theta_n}V          @VV{1}V \\
\K_{n,1}  @>{j}>> \Emb_*(S^1,S^n)  @>{}>> S^{n-1}.
\end{CD}
\]

The map $ \Theta_n: SO(n) \times_{SO(n-1)} \K_{n,1} \to\
\Emb_*(S^1,S^n)$ is thus a homotopy equivalence. Hence the map
$$\Theta_n: SO(n+1) \times_{SO(n-1)} \K_{n,1}\to\ \Emb(S^1,S^n)$$
is also a homotopy equivalence.
\end{proof}

Restrict attention to the special case given by $n=3$.

\begin{cor}\label{neq3long_closed_rel}
There is a homeomorphism
$$S^3 \times \Emb_*(S^1,S^3) \to \Emb(S^1,S^3).$$ Furthermore, the bundle
$$\K_{3,1} \to \Emb_*(S^1,S^3) \to S^2$$
is the induced bundle with fibre $\K_{3,1}$ from the bundle
$$SO(2) \to SO(3) \to S^2$$ where $SO(2)$ acts on $\K_{3,1}$ by rotation
about the long axis, as previously described. Thus, up to a
homotopy-equivalence $\Emb_*(S^1,S^3)$ is the union of two copies of
$D^2 \times \K_{3,1}$ along their common boundary, where the gluing
map $S^1 \times \K_{3,1} \to S^1 \times \K_{3,1}$ is given by
$(z,f)\longmapsto (z,z^2.f) \in S^1 \times \K_{3,1}$, where we
identify $S^1 \equiv SO(2)$ and its action by rotation about the
long axis.
\end{cor}

Observe that Proposition \ref{closed-long-prop} generalizes to a
proposition about the embedding spaces $\Emb(S^k,S^n)$. 
We skip the proof as it is essentially the same as Proposition
\ref{closed-long-prop}.

\begin{prop}\label{skinsn}
Provided $n-j \geq 1$, there is a homotopy-equivalence
$SO(n+1)\times_{SO(n-j)} \K_{n,j} \to \Emb(S^j,S^n)$.
\end{prop}

Given $f \in \K_{n,j}$, let $\dot f \in \Emb(S^j,S^n)$ be the
one-point compactification of $f$. Consider $S^n$ to be the one-point
compactification of $\Real^n$. The inclusion $\Real^n \to S^n$
induces an inclusion $\Emb(S^j,\Real^n) \to \Emb(S^j,S^n)$.
 This inclusion induces a bijection $\pi_0 \Emb(S^j,\Real^n) \to \pi_0
\Emb(S^j,S^n)$ provided $n-j \geq 2$. Given $f \in \K_{n,j}$ let
$\bar f \in \Emb(S^j,\Real^n)$ be such that $\dot f$ is isotopic (in
$S^n$) to $\bar f$. These conventions give us a one-to-one
correspondence between $\pi_0 \K_{n,j}$, $\pi_0 \Emb(S^j,S^n)$ and
$\pi_0 \Emb(S^j,\Real^n)$ for $n-j \geq 2$.

If $f \in \K_{n,j}$ is a long knot, let $X_f$ denote the component
of $\bar f$ in $\Emb(S^j, D^n)$ and let $C_f$ denote the complement
of an open tubular neighbourhood of $\dot f$ in $S^n$. 
Given $f \in \K_{n,j}$ define
$C_f \rtimes \K_{n,j}(f) =
 \{(p,g) : g \in \K_{n,j}, p \in C_g, \text{ where } g \text{ isotopic to } f\}$, 
and define $C \rtimes \K_{n,j}$ to be the union of the spaces $C_f \rtimes \K_{n,j}(f)$
for all $f \in \K_{n,j}$

\begin{prop}\label{long-cl} Provided $n-j > 0$, $\Emb(S^j,\Real^n)$ is homotopy-equivalent to the
space $SO(n) \times_{SO(n-j)} \left( C \rtimes \K_{n,j} \right)$.  In particular the
components $X_f$ of $\Emb(S^j,\Real^n)$ have the homotopy-type of
$SO(n) \times_{SO(n-j)} \left( C_f \rtimes \K_{n,j}(f) \right)$.
$SO(n-j) \subset SO(n)$ acts on $SO(n)$ as the subgroup fixing
a $j$-dimensional subspace, and $SO(n-j)$ acts on $C_f \rtimes
\K_{n,j}$ diagonally.
\begin{proof}
See Proposition 2.2 of \cite{B3}.
\end{proof}
\end{prop}

\section{On the homology of $\Cu_2(X \amalg \{*\})$}\label{c2homol}

The purpose of this section is to recall the homology of $$\Cu_2(X
\amalg \{*\})$$ for $X$ not necessarily path-connected. These
results will then be combined with Theorem \ref{long knots in 3
space} to obtain Theorem \ref{homology of long knots in 3 space}.
The space $X$ is assumed to be compactly generated and weak
Hausdorff as a topological space \cite{CLM}; the base-point $\{*\}$
is non-degenerate by construction.

Formal constructions are given next for which $\F$ is a field and
all modules are assumed to be vector spaces over $\F$. Let $V$
denote a graded vector space which splits as a direct sum $$V = V_+
\oplus V_{-}$$ for which $V_+$ consists of the elements concentrated
in even degrees and $V_{-}$ consists of the elements concentrated in
odd degrees. Let $\sigma(V)$ denote the ``algebraic suspension" of
$V$. That is $\sigma(V)$ is the module $V$ with all degrees raised
by one. In addition, define $\sigma^n(V) =
\sigma(\sigma^{n-1}(V)).$ The ``algebraic desuspension" of $V$
denoted $\sigma^{-1}(V)$ is defined by requiring
$\sigma(\sigma^{-1}(V)) = V.$

Next consider the free Lie algebra $$L[\sigma(V)]$$ and, if $\F =
\F_p$, the free restricted Lie algebra over $\F_p$ denoted
$L^{(p)}[\sigma(V)].$ In this last case, consider the natural
inclusion $$j:L[\sigma(V)] \to L^{(p)}[\sigma(V)]$$ with co-kernel
denoted $W^p[\sigma(V)]$ (for which $L[\sigma(V)]$ is the free Lie
algebra defined over the field $\F_p$).  The definition of a
restricted Lie algebra is given in Jacobson's book ``Lie Algebras"
\cite{Jac} with graded restricted Lie algebras treated in \cite{MM}.
Graded restricted Lie algebras may be regarded as the module of
primitive elements in the tensor algebra $T[\sigma(V)]$ defined over
the field $\F_p$. Notice that $\sigma^{-1}(\sigma V) = V$,
but that $\sigma^{-1}L[\sigma(V)]$ is not isomorphic to $L[V]$ in
general.

Let $E[V_-]$ denote the exterior algebra generated by $V_-$ and let
$\F[V_+]$ denote the polynomial algebra generated by $V_+$. Consider
the symmetric algebra $S[V]$ defined as follows:
\begin{enumerate}
    \item For the field $\mathbb Q$, $S[V]= E[V_{-}] \otimes \F[V_{+}]$.
    \item For the field $\F_2$, $S[V]= \F[V]$, the polynomial
    algebra generated by $V$.
    \item For the field $\F_p$ with $p$ an odd prime, 
    $S[V]= E[V_{-}] \otimes \F[V_{+}]$.
\end{enumerate}

We describe the homology of $\Cu_2(X \amalg \{*\})$ with coefficients in
the field (1) $\mathbb Q$, (2) $\F_2$ and (3) $\F_p$ for $p$ and odd prime.

Consider case (1). Let $V= H_*(X,\mathbb Q),$ and form the symmetric algebra 
$S[\sigma^{-1}L[\sigma(H_*(X; \mathbb Q)]].$ By \cite{CLM}, there is
an isomorphism of Hopf algebras $$S[\sigma^{-1}L[\sigma(H_*(X;
\mathbb Q))]] \to\ H_*(\Cu_2(X \amalg \{*\}); \mathbb Q)$$ with
co-product determined by that of $H_*(X; \mathbb Q)$.

The analogous theorem for $\F_2$ is given as follows. Let
$V= H_*(X;\F_2 ),$ and form the symmetric algebra
$$ S[\sigma^{-1}L^{(2)}[\sigma(H_*(X; \F_2))]].$$ By \cite{CLM},
there is an isomorphism of Hopf algebras
$$S[\sigma^{-1}L^{(2)}[\sigma(H_*(X;\F_2))]] \to\ H_*(\Cu_2(X \amalg
\{*\}); \F_2 )$$ with co-product determined by that of
$H_*(X;\F_2)$. Remark: The role of the restriction in a restricted
Lie algebra over $\F_2$ is to create the Araki-Kudo-Dyer-Lashof
operation, the operation which sends an element $\sigma(v)$ to
$Q_1(\sigma(v))$.

The result for odd primes $p$ with $\F_p$ is given as
follows. Let $V= H_*(X;\F_p ),$ and form the symmetric algebra
$$S[\sigma^{-1}L^{(p)}[\sigma(V)] \oplus \sigma^{-2}W^p[\sigma(V)]].$$
By \cite{CLM}, there is an isomorphism of Hopf algebras
$$S[\sigma^{-1}L^{(p)}[\sigma(V)] \oplus \sigma^{-2}W^p[\sigma(V)]] \to\ H_*(\Cu_2(X \amalg
\{*\}); \F_p )$$ with coproduct determined by that of $H_*(X;\F_p)$.

\section{On the homology of $\K_{3,1}$}\label{khomol}

Recall the  homotopy equivalence of Theorem \ref{long knots in 3
space},

$$C(\Real^2,X_{\K} \amalg \{*\} ) \to \K_{3,1}$$

Thus there are isomorphisms of Hopf algebras
\begin{enumerate}
    \item $S[\sigma^{-1}L[\sigma(V)]] \to H_*(\K_{3,1}; \mathbb Q)$
    for $V = H_*(X_{\K}; \mathbb Q)$,
    \item $S[\sigma^{-1}L^2[\sigma(V)]]\to H_*(\K_{3,1}; \F_2)$ for
    $V = H_*(X_{\K};\F_2)$ and
    \item $S[\sigma^{-1}L^{(p)}[\sigma(V)] \oplus \sigma^{-2}W^p[\sigma(V_{-})]]$
for $V = H_*(X_{\K};\F_p)$ for odd primes $p$.
\end{enumerate}

Further information  concerning the the homology of $X_{\K}$ is
given in Section \ref{tcns}.

Thus the above isomorphisms give the homology of $\K_{3,1}$ with
field coefficients. The first and second parts of Theorem
\ref{homology of long knots in 3 space} follow. The proof of Theorem
\ref{homology of long knots in 3 space} will be concluded in Section
\ref{findtorsion} in which higher torsion is constructed.

Notice that the space $\Cu_2(X \amalg \{*\})$ is naturally a
disjoint union of $X$ with another space. Thus, there is a natural
direct sum decomposition of graded vectors spaces
$$\bar H_*(X; \mathbb F) \oplus \Gamma(X; \F) \to\
 H_*(\Cu_2(X \amalg \{*\}); \mathbb F))$$ for a choice
of graded vector space $$\Gamma(X; \F)$$ which is functor of $H_*(X;
\F)$.

The construction $\Gamma(X; \F)$ is used in Section \ref{tcns} to
describe the homology of the subspace of $\K_{3,1}$ generated by
torus knots, as well as the operations of connected sums, cablings
and the action of the little two-cubes.

\section{Higher $p$-torsion in the integer homology of $\K_{3,1}$}\label{findtorsion}

One way in which higher order $p$-torsion in the homology of
$\K_{3,1}$ arises is summarized next. The way in which little cubes
$\Cu_2(n)$ are related to configuration spaces $\cfg(\Real^2,n)$ is
as follows. There are maps $\Cu_2(n) \to \cfg(\Real^2, n)$ which are
both homotopy equivalences and equivariant with respect to the
action of the symmetric group $\Sigma_n$ \cite{M}. Thus it suffices
to exhibit higher torsion in the integer homology of
$\cfg(\Real^2,n) \times_{\Sigma_{n}}X^{n}$ for certain choices of
spaces $X = \K_{3,1}({f})$. Since the construction $\cfg(\Real^2,n)
\times_{\Sigma_{n}}X^{n}$ occurs numerous times below, it is
convenient to define $$E_n(X) = \cfg(\Real^2,n)
\times_{\Sigma_{n}}X^{n}$$ as given in the Introduction.

Given a prime long knot $f$, consider the path-component
$\K_{3,1}(g)$ where
$$g = \#_{n}f \equiv f \rep_n \splice \Keychain^n$$
Here $\#_{n}f$ denotes the connected-sum of $n$ copies of the same
knot $f$, with the splice notation from Section \ref{ops} given.
There are homotopy equivalences $$\K_{3,1}(g) \to \cfg(\Real^2, n)
\times_{\Sigma_{n}} {\K_{3,1}(f)}^{n} = E_L(\K_{3,1}(f)).$$

First consider $p$-torsion of order exactly $p$ obtained from the
equivariant cohomology of $\cfg(\Real^2, p)$ as constructed in
\cite{CLM}.

\begin{prop}\label{summand}
Let $Y$ denote any connected CW-complex.
\begin{enumerate}
    \item If $H_{2t-1}(Y; \F_p)$ is non-zero, then $\F_p$
is a direct summand of $$H_{2pt-2}(\cfg(\Real^2, p)
\times_{\Sigma_{p}} {Y}^{p}; \Zed).$$
    \item If $\cyc{p^s}$ is a direct summand of
    $H_{2t-1}(Y; \Zed )$, then
    $$H_{2tp^{r}-1}(E_{p^{r}}(Y);\Zed) = H_{2tp^{r}-1}(\cfg(\Real^2,p^{r})
    \times_{\Sigma_{p^{r}}}Y^{p^{r}}; \Zed)$$
has a $\cyc{p^{s+q}}$-summand.

 \item There is a homotopy equivalence $$E_n(\K_{3,1}(f)) \to
\K_{3,1}(\#_{n}f).$$

   \item Thus if $\K_{3,1}(f)$ has any
non-trivial mod-$p$ homology in degree $2t-1$, then $\F_p$ is a
direct summand of $$H_{2pt-2}(E_p(\K_{3,1}(f)); \Zed).$$ Hence
$$H_{2pt-2}(\K_{3,1}(\#_{p}f);\Zed) = \F_p \oplus A$$ for some
abelian group $A$.

\item Furthermore, if $\cyc{p^s}$ is a direct
summand of $H_{2t-1}(\K_{3,1}(f); \Zed )$, then
$$H_{2tp^{r}-1}(E_{p^{r}}(\K_{3,1}(f));\Zed) = H_{2tp^{r}-1}(\K_{3,1}(\#_{p^r}f);\Zed)
= \cyc{p^{s+r}} \oplus A$$ for some abelian group $A$.

\end{enumerate}

\end{prop}

Assume that $$H_{j}( \K_{3,1}(f); \Zed) = \cyc{p^s} \oplus A$$ for
some abelian group $A$. Label this $\cyc{p^s}$-summand
(ambiguously) by $<f,j,\cyc{p^s}>.$

\begin{exm}\label{torus knot example}
By \cite{Hatcher4} or \cite{B2} if $f$ is a non-trivial torus knot,
then $\K_{3,1}(f)$ has the homotopy type of a circle. A direct
application of Proposition \ref{summand} gives that
$$H_{2p-2}(\cfg(\Real^2, p) \times_{\Sigma_{p}}
{\left(\K_{3,1}(f)\right)}^{p}; \Zed) = H_{2p-2}( \K_{3,1}(\#_{p}f);
\Zed) = <\#_{p}f,2p-2,\F_p> \oplus A$$ for some abelian group $A$.
\end{exm}

Next, consider the prime knot $h$ given by a cabling of $f$.
\begin{exm}\label{iterated cabling example} The
examples here arise by an iterated cabling construction as follows.
Let $\alpha/\beta \in \Rational$ satisfy $\beta \geq 1$ with
$GCD(\alpha,\beta) = 1$. $h=f \splice \Seifert^{(\alpha,\beta)}$ is
the $\alpha/\beta$-cabling of $f$.

There is a homotopy equivalence $$\K_{3,1}(f)\times S^1
\to\K_{3,1}(h) = \K_{3,1}(f \splice \Seifert^{(\alpha,\beta)})  $$
as described in Section \ref{htpe}.

Next, consider an $m$-fold iterated cabling $h_m$ of $f$ defined by
\begin{itemize}
    \item $h_1 = f \splice \Seifert^{(\alpha,\beta)}$,
    \item $h_{i+1} = h_i \splice \Seifert^{(\alpha,\beta)}$ for $i \geq 1$, defined recursively.
\end{itemize}

Then there are homotopy equivalences $$\K_{3,1}(h_m) \to \K_{3,1}(f)
\times (S^1)^m.$$

Assume that $$H_{2t-1}(\K_{3,1}(f));\mathbb F_p)= <f,2t-1,\F_p>
\oplus A$$ for some abelian group $A$. Then the integer homology of
$\K_{3,1}(h_m) = \K_{3,1}(f) \times (S^1)^m$  has a summand denoted
(ambiguously) by
$$<f,2t-1,\F_p> \otimes H_*( (S^1)^m;\Zed).$$

Thus if $ m \geq 1$, there are elements of order $p$ in both odd as
well as even degrees in the integer homology of $\K_{3,1}(h_m))$. If
$m$ is ``large", then there are many elements of order exactly $p$
which are of both odd and even degree. This fussiness concerning
parity of degrees has consequences for higher torsion in homology.
\end{exm}

The above remarks give examples of long knots with torsion of order
exactly $p$ concentrated in odd degrees for the integer homology of
their path-components. One choice of $f$ is a torus knot. The next
proposition shows that $p$-torsion of order exactly $p$ in the
homology of $\K_{3,1}(g)$ gives rise higher $p^s$-torsion in the
homology of other components related components as follows.

Recall that Example \ref{iterated cabling example} provides
instances of prime knots $f$ such that
$$H_{2t-1}(\K_{3,1}(f);\Zed) = <f,2t-1,\F_p>
\oplus A$$ for some abelian group $A$. Consider the long knot $g$
given by the $p^s$-fold connected sum
$$g = f \rep_{p^s} \splice \Keychain^{p^s} \equiv \#_{p^s}f $$ as used
in the next Proposition.
\begin{prop}\label{higher p-tosion example}
Let $f$ denote a prime long knot such that
$$H_{2t-1}(\K_{3,1}(f);\Zed) = <f,2t-1,\F_p>
\oplus A$$ for some abelian group $A$. Let
$$g = f \# \cdots \# f = \#_{p^s}f.$$
Then $$H_{2tp^s -1}(\K_{3,1}(g);\Zed) = <\#_{p^s}f,2tp^s
-1,\cyc{p^{s+1}}> \oplus A$$ for some abelian group $A$.
\end{prop}

\begin{proof}

Assume that the integer homology of $ \K_{3,1}(f)$ has a non-trivial
$\F_p$ summand in degree $2t-1$ as guaranteed by example
\ref{iterated cabling example}. Thus, in the mod-$p$ reduction of
the integer homology of $\K_{3,1}(f)$, there are classes $x$ of
degree $2t-1$, the mod-$p$ reduction of a class of order $p$, as
well as a class $y$ in degree $2t$ which corresponds to the
contribution forced by $x$ in the ``Tor" term in the classical
universal coefficient Theorem.

Since $H_{2t-1}(\K_{3,1}(f);\Zed)$ has a $\F_p$-summand, there are
elements in the mod-$p$ homology of $\K_{3,1}(f)$ with

\begin{enumerate}
    \item $x$ in $H_{2t-1}(\K_{3,1}(f);\F_p)$,
    \item $y$ in $H_{2t-1}(\K_{3,1}(f);\F_p)$ and
    \item the first Bockstein of $y$ is $x$, $$\beta_{1}(y) = x.$$
\end{enumerate}

A classical computation of the Bockstein spectral sequence gives
that the $(s+1)$-st Bockstein is defined with
$$\beta_{s+1}(y^{p^s}) = x \cdot y^{-1 + p^s} + I
$$ in case $s\geq 1$ for which $I$ denotes the indeterminacy of this
operation. The Proposition follows as these classes survive in the
Bockstein spectral sequence for $C(\Real^2,X_{\K} \amalg \{*\} )$.
\end{proof}

Two concrete examples are listed next.
\begin{exm}\label{higher torsion}

    \begin{enumerate}
    \item Let $f$ denote a non-trivial torus knot with
    $$H_{1}\K_{3,1}(f) =\Zed.$$ Then the long knot
$$\#_{p}f = f \# \cdots \# f $$ satisfies the property
that $$ \K_{3,1}(\#_{p}f) = E_p(\K_{3,1}(f))$$ with
$$H_{2p-2}(\K_{3,1}(\#_{p}f);\Zed) = <\#_{p}f,2p-2,\cyc{p^1}> \oplus A$$ for some abelian group $A$.

\item Denote a cable of $\#_{p}f$
by $(\#_{p}f) \splice \Seifert^{\alpha,\beta}$. There are homotopy
equivalences
$$\K_{3,1}((\#_{p}f)\splice \Seifert^{\alpha,\beta}) \to \K_{3,1}(\#_{p}f) \times S^1$$
with the property that
$$H_{2p-1}(\K_{3,1}((\#_{p}f)\splice \Seifert^{\alpha,\beta} );\Zed)
= < (\#_{p}f) \splice \Seifert^{\alpha,\beta},2p-1, \F_p> \oplus
A$$ for some abelian group $A$.

\item The $p^s$-fold connected sum of $(\#_{p}f)\splice \Seifert^{\alpha,\beta}$,
$$\#_{p^s} \left( (\#_{p}f)\splice \Seifert^{\alpha,\beta}\right) $$ has the property
that
$$H_{2p^{s+1}-1}\left(\K_{3,1}\left(\#_{p^s} \left( (\#_{p}f)\splice \Seifert^{\alpha,\beta}\right) \right);\Zed\right) = <\#_{p^s} \left( (\#_{p}f) \splice
\Seifert^{\alpha,\beta}\right) ,2p^{s+1}-1,\cyc{p^{s+1}}> \oplus
A$$ for some abelian group $A$.
\end{enumerate}
\end{exm}

The third part of Theorem \ref{homology of long knots in 3 space}
follows, thus concluding the proof.

\section{On $H_1 \K_{3,1}$}
\label{initial_computations_of_H}

$H_* \K_{3,1}$ has torsion of all orders, it is natural to ask for
the lowest dimension $i_{(p,n)}$ so that $H_{i_{(p,n)}} \K_{3,1}$
contains torsion of order $p^n$.  This question is answered in this
section for the special case $(p,n)=(2,1)$.  This section contains a
proof of Theorem \ref{H1.of.the knot.space}.

The idea of the proof is to describe $H_1 \K_{3,1}$ inductively,
component-by-component.  The most complicated case from the point of
view of torsion is the hyperbolic satellite case, since there is
currently insufficient control of the representations $B_L \to
\Sigma_n^+$. In addition, better control over the class of the
inversion-map $H_1 \K_{3,1} \to H_1 \K_{3,1}$ is required.

First, the base-case: knots whose JSJ-trees have one vertex.

\begin{itemize}
\item If $f$ is a torus knot $H_1 \K_{3,1}(f) \simeq \Zed$ and $f$ is invertible
and the $\Sigma_2$ action (inversion action on $H_1 \K_{3,1}(f)$) is given by
multiplication by $(-1)$.  This is a direct corollary of \cite{B2}.

\item If $f$ is a hyperbolic knot $H_1 \K_{3,1}(f) \simeq \Zed^2$.  In the
case that $f$ is invertible, the (inversion) action of $\Sigma_2$ on
$\Zed^2$ is multiplication by $(-1)$. This also follows immediately
from \cite{B2}.
\end{itemize}

The next proposition gives $H_1 (\K_{3,1}(f);\Zed)$ inductively, via
the JSJ-tree of $f$. Given a group $G$ acting on an abelian group
$A$, let $A_G$ denote the module of co-invariants, the quotient
group of $A$ modulo the subgroup generated by $\{a-g\cdot a : g \in
G, a \in A\}$. The following lemma follows from the Leray-Serre
spectral sequence for any fibration with a section.

\begin{lem}\label{lem_h1} Given a fibration $F \to E \to B$ with
a section, with both the base and the fibre path-connected, then
$H_1 E \simeq H_1 B \oplus \left(H_1 F\right)_{\pi_1 B}$.
\end{lem}

In principle, one can deduce the following result
from the presentation of the groups $\pi_1 (\cfg(\Real^2,n)/Y;\Zed)$
given by Manfredini \cite{Manfredini}, alternatively using some
elementary facts about the abelianization of the braid group or from
the description in \cite{CLM}. The proof is omitted.

\begin{lem}\label{lem_ab}
Let $Y$ be a Young subgroup of $\Sigma_n$. We think of $\Sigma_n$ as the group
of bijections of the set $\{1,2,\cdots,n\}$. Let $k$ be the number
of distinct orbits of $Y$, and let $f$ be the number of fixed-points
of $Y$ acting on $\{1,2,\cdots,n\}$. Let $l=k-f$. Then $H_1
(\cfg(\Real^2,n)/Y;\Zed)$ is a free-abelian group of rank $l+{k
\choose 2}$.
\end{lem}

\begin{prop}
Given any component $\K_{3,1}(f)$ of $\K_{3,1}$, $H_1
(\K_{3,1}(f);\Zed)$ is finitely-generated and a direct-sum of groups
of the form: $\Zed$ and $\cyc{2}$. Moreover, if $f$ is an invertible
knot, the involution of $H_1 (\K_{3,1}(f);\Zed)$ 
preserves a splitting of $H_1 (\K_{3,1}(f);\Zed)$ into a
direct sum $H_1 (\K_{3,1}(f);\Zed) = V_1 \oplus V_2$ where the
involution acts on $V_1$ as the identity and acts on $V_2$ by
multiplication by $(-1)$.
\begin{proof}

The proof is by induction on the height of the JSJ-tree of $f$. The
height one case was dealt with at the start of this section.
The inductive step is as follows.

\begin{itemize}
\item Consider the cases that $f$ is a cable of $g$, then
$H_1 \K_{3,1}(f) = \Zed \oplus H_1 \K_{3,1}(g)$ \cite{B2}.  In the
case that $f$ is invertible, the homotopy-equivalence is
$\F_2$-equivariant with $\F_2$-action on $S^1 \times \K_g$ being
complex conjugation on $S^1$ and the inversion involution on $\K_g$.
Thus the $\F_2$-action on $H_1 S^1 \times \K_{3,1}(g) = \Zed
\oplus H_1 \K_{3,1}(g)$ is simply the direct sum of $\F_2$-modules
$\Zed$ (with the non-trivial involution) and $H_1 \K_{3,1}(g)$ with
its own inversion involution, completing this part of the inductive
step.

\item Consider the case that $f$ is a connected sum of prime knots
$f_1, f_2, \cdots, f_n$ with $n \geq 2$.  Then by Lemma
\ref{lem_h1},

$$H_1 (\K_{3,1}(f);\Zed) = H_1(\cfg(\Real^2,n)/Y;\Zed) \oplus
\left( \bigoplus_{i=1}^n H_1 (\K_{3,1}(f_i);\Zed) \right)/Y$$ where
$Y$ is the Young subgroup of $\Sigma_n$ given by the equivalence relation
$i \sim j \Leftrightarrow \K_{3,1}(f_i)=\K_{3,1}(f_j)$. $H_1
(\cfg(\Real^2,n)/Y;\Zed) \simeq \Zed^{l + {k \choose 2}}$ where $l$
is the number of orbits of $Y$ with more than $1$ element, and $k$
is the number of orbits of $Y$ by Lemma \ref{lem_ab}. If $f$ is
invertible, the involution action on $\K_{3,1}(f) \simeq
\cfg(\Real^2,n) \times_{Y} \prod_{i=1}^n \K_{3,1}(f_i)$ was
described in \cite{B2} as a map that preserved the above product
structure, acting by mirror reflection along a line in $\Real^2$ on
$\cfg(\Real^2,n)$ and by permutation of the factors of
$\K_{3,1}(f_i)$ combined with the inversion involution on
$\K_{3,1}(f_i)$ for each $i \in \{1,2,\cdots,n\}$. Since the
abelianisation of $\pi_1 (\cfg(\Real^2,n)/Y)$ was computed entirely
in terms of linking numbers, mirror reflection along a line induces
multiplication by $(-1)$ on $H_1 (\cfg(\Real^2,n)/Y;\Zed)$. This
completes this step of the inductive argument.

\item Consider the case of a hyperbolic satellite operation.  In this case,
$H_1 \K_f = \Zed^2 \oplus \left( \oplus H_1 \K_{3,1}(f_i) \right) /
A_f$ \cite{B2}.  Thus $H_1 \K_f$ consists of $\Zed^2$ direct sum
various groups, one for each orbit of $A_f$ acting on
$\{1,2,\cdots,n\}$. Denote the orbits by $\{1,2,\cdots,n\}=Y_1 \cup
Y_2 \cup \cdots \cup Y_k$. The summand corresponding to orbit $Y_i$
is either $H_1 \K_{3,1}(f_j)$ for $j \in Y_i$ or $(H_1
\K_{3,1}(f_j))/\Sigma_2$ depending on whether or not $A_f$ has an
element that reverses the orientation of $L_j$ or not.  $(H_1
\K_{3,1}(f_j))/\Sigma_2$ is also a direct sum of groups of the form
$\Zed$ and $\cyc{2}$ by the inductive hypothesis. Now consider the
case that $f$ is invertible. By \cite{B2} the $\Sigma_2$-action on
$\K_f \simeq M \times (SO(2) \times_{A_f} \prod_{i=1}^n
\K_{3,1}(f_i))$ respects the bundle structure, thus on $\Zed^2
\oplus \left( \oplus H_1 \K_{3,1}(f_i) \right) / A_f$ it acts by
multiplication by $(-1)$ on the $\Zed^2$-factor.  On the remaining
factors it either acts trivially on the $i$-th summand if the
inversion symmetry of $L$ does not reverse the orientation of $L_i$
or it acts by inversion on that summand.
\end{itemize}
Thus the result follows.
\end{proof}
\end{prop}

\begin{cor}
$H_1 (\K_{3,1}(f);\Zed)$ contains $2$-torsion if and only if there
is a hyperbolic link $L$ so that one of the vertices of $\IG_f$ is
decorated by $L$, and if we let $g$ be the knot whose JSJ-tree is
the subtree rooted at $L$, then $A_g$ must contain an isometry that
reverses the orientation of some $L_i$.
\end{cor}

\section{The first occurrence of odd torsion}
\label{no_p_torsion}

\begin{thm}\label{first_oddp_torsion}
Let $f$ denote a long knot and $p$ an odd prime. If
$H_i(\K_{3,1}(f);\Zed)$ contains $\F_p$, then $i \geq 2p-2$.
\end{thm}

Notice that the Theorem does not assert that there is torsion in
$H_* \K_{3,1}(f)$, but rather the least dimension in which
$p$-torsion can possibly occur. There are long knots $f$ such that
$H_* \K_{3,1}(f)$ is torsion free. Furthermore, there are long knots
$g$ such that $H_{2p-2}(\K_{3,1}(g);\Zed)$ contains copies of
$\F_p$ by Example \ref{torus knot example}. This theorem follows a
classical pattern which is exhibited for both $K(B_n,1)$ as well as
$\Omega^n S^{n+1}$.

\begin{proof}
It suffices to prove that $H_i(\K_{3,1};\Zed)$ contains no
$p$-torsion for $i<2p-2$.

Since the homology groups $H_* \K_{3,1}(f)$ are torsion free for the
unknot, torus knots and hyperbolic knots, it suffices to check that
that $p$-torsion cannot occur in dimensions less than $2p-2$ in the
following three cases.

\begin{enumerate}
  \item The knot $f$ is a cable of $g$ in which case $\K_{3,1}(f) \simeq S^1 \times
\K_{3,1}(g)$.
  \item The knot $f$ is hyperbolically spliced.
  \item The knot $f$ is a connected-sum of knots $g_i$ such that the
  homology of $\K_{3,1}(g_i)$ is $p$-torsion free in dimensions
  less than $2p-2$.
\end{enumerate}
Case $1$ follows directly from the classical K\"unneth theorem.
Cases $2$ and $3$ follow inductively by the next three lemmas.
\end{proof}

Consider the $k$-fold product $X^k$ with the natural (left) action
of $\Sigma_k$ on $X^k$.  The free $\Cu_n$-space
generated by $X\amalg +$ is denoted $\Cu_n(X\amalg +)$. Recall that
this space is homotopy-equivalent to the disjoint union of 
$\cfg(\mathbb R^n,k) \times_{\Sigma_k}X^k$ for all $k \geq 0$.

\begin{lem}\label{lem:torsion.for.covering spaces}
Assume that $X$ is a topological space of the homotopy type of a
CW-complex (alternatively, one can substitute compactly generated,
weak Hausdorff for having the homotopy-type of a CW-complex in this
lemma) without $p$-torsion in homology of dimensions less than
$2p-2$ for $p$ an odd prime. Then the homology of $\Cu_n(X\amalg +)$
and $\Omega^n \Sigma^n(X\amalg +)$ also do not have $p$-torsion in
homology of dimensions less than $2p-2$.
 Thus the homology of $\cfg(\mathbb
R^n,k) \times_{\Sigma_k}X^k$ does not have $p$-torsion in homology
of dimensions less than $2p-2$.
\end{lem}

The proof follows directly from the computations in \cite{CLM} or
can be done classically by chain level arguments (in the spirit of
Nakaoka and Steenrod).

\begin{lem}\label{lem:bundles.over.circles}
Let $A = \mathbb Z/p^r\mathbb Z$ act on the $p^r$-fold product of a
path-connected CW-complex $X^{p^r}$ by a cyclic permutation of order
$p^r$ and on $S^1$ freely via a rotation of order $p^r$. If
$H_i(X;\mathbb Z)$ is $p$-torsion free and finitely generated for
all $i< q$, then $$H_j( S^1 \times_{A}X^{p^r}; \mathbb Z)$$ is
$p$-torsion free for all $j < q$.
\end{lem}

\begin{proof}
Consider the space $$S^1 \times_{A}X^{p^r}.$$  Classically, there
are chain equivalences

\[
\begin{CD}
\mathbb B_* \otimes_{\mathbb Z[A]} C_*(X)^{\otimes p^r} @>{D \otimes
1^{p^r}}>> \mathbb B_* \otimes_{\mathbb Z[A]} C_*(X^{p^r}) @>{} >>
C_*(S^1 \times_{A}X^{p^r}  )
\end{CD}
\] where

\begin{enumerate}
\item $A = \mathbb Z/p^r \mathbb Z$,
\item $\mathbb B_*$ denotes the chain complex of right
$\mathbb Z[A]$-modules ( chain equivalent to the total singular
chain complex of a circle )

\[
\begin{CD}
\cdots @>{} >> \{0\} @>{} >> \mathbb Z[A] @>{D}>> \mathbb Z[A]
\end{CD}
\] for which $D(1) = 1-\tau$ where $\tau$ is a generator for $A$,

\item $C_*(X)$ denotes the total singular chain complex of $X$ for
which $C_i(X)$ denotes the singular chains in degree $i$ and
\item $C_*(X)^{\otimes p^r}$ and $C_*(X^{ p^r})$ is given the natural
structure of left $\mathbb Z[A]$-modules.
\end{enumerate}

Since $X$ is assumed to be of finite type, $$H_i(X) = F_i \oplus
T_i$$ where $F_i$ is a finite direct sum of copies of $\mathbb Z$
and $T_i$ is a finite direct sum of finite cyclic groups. If $ i<
q$, it may be assumed that $T_i$ is of order prime to $p$ and thus
this summand will not contribute $p$-torsion to the homology of the
chain complex $\mathbb B_* \otimes_{\mathbb Z[A]} C_*(X)^{\otimes
p^r}$ ( Details are omitted ).

Furthermore, there is a map $$\rho_i: F_i \to C_i(X)$$ which
\begin{enumerate}
  \item induces a map of chains complexes ( with trivial differential
  for the source ) and
  \item induces a homology isomorphism in degrees $i < q$ with coefficients
  in $\mathbb Z_{(p)}$ ( the integers localized at $p$ meaning those
  rational numbers with denominators prime to $p$ ).
\end{enumerate}

Thus it suffices to check that the homology of the chain complex
$$ \mathbb B_* \otimes_{\mathbb Z[A]} (\oplus_{i < q}
F_i)^{\otimes p^r}$$ is $p$-torsion free homology in dimensions less
than $q$ for $p$ an odd prime. Since $F_i$ is free abelian, notice
that $(\oplus_{i < q} F_i)^{\otimes p^r}$ is a sum of permutation
representations ( over the integral group ring of $A$ ) each of
which are cyclic $\mathbb Z[A]$-modules which have the following
generators.

\begin{enumerate}
  \item $v^{\otimes p^r}$ where $v$ is an element in $F_i$ of even
  degree.
  \item $v^{\otimes p^r}$ where $v$ is an element in $F_i$ of odd
  degree.
  \item $v_1 \otimes v_2 \otimes \cdots \otimes v_{p^r}$ where the
  $v_i$ run over a basis for the $\oplus_{i < q}F_i$ with at least two distinct
  basis elements appearing.
\end{enumerate}

Thus it suffices to work out the torsion in the chain complex
$$ \mathbb B_* \otimes_{\mathbb Z[A]} M $$ where $M$ denotes the
free abelian group which is a cyclic $\mathbb Z[A]$-module with one
of the elements in (1-3) as generators. These are considered next.

\begin{enumerate}

\item Let $M$ denote the cyclic $\mathbb Z[A]$-module generated
by $v ^{\otimes p^r}$ where $v$ is an element in $F_i$ of either odd
or even degree. Since $p$ is odd, the associated permutation
representation is trivial and thus the chain complex
$$\mathbb B_* \otimes_{\mathbb Z[A]} M$$ is isomorphic to $(\mathbb
B_* \otimes_{\mathbb Z[A]} \mathbb Z) \otimes_{\mathbb Z} M$. The
homology of this chain complex is isomorphic to $H_*(S^1)
\otimes_{\mathbb Z} M$ as a graded abelian group and is thus torsion
free.

\item Let $M$ denote the cyclic $\mathbb Z[A]$-module generated
  by  $v_1 \otimes v_2 \otimes \cdots \otimes v_{p^r}$ where the
  $v_i$ run over a basis for the $\oplus_{i < q}F_i$ with at least two distinct
  basis elements appearing among the $v_i$. The action of $A = \mathbb Z/p^r \mathbb
  Z$ has isotropy subgroup given by $\mathbb Z/p^s \mathbb
  Z$ for some $ 0 \leq s < r$. Thus the module $M$ is isomorphic to
  to $$\mathbb Z[A] \otimes_{\mathbb Z[\mathbb Z/p^s \mathbb Z] }\mathbb Z$$
as a left $\mathbb Z[A]$-module and there is an induced isomorphism
of chain complexes
$$\mathbb B_* \otimes_{\mathbb Z[A]} (\mathbb Z[A] \otimes _{\mathbb Z[\mathbb Z/p^s \mathbb Z] }\mathbb
Z)\to \mathbb B_* \otimes_{\mathbb Z[A]} M.$$

Since the chain complex $\mathbb B_* \otimes_{\mathbb Z[A]} (\mathbb
Z[A] \otimes _{\mathbb Z[\mathbb Z/p^s \mathbb Z] }\mathbb Z)$ is
isomorphic to $$\mathbb B_* \otimes _{\mathbb Z[\mathbb Z/p^s
\mathbb Z] }\mathbb Z,$$ the chain complex has torsion free
homology.
\end{enumerate}

\end{proof}

\begin{lem}\label{lem:birth.of.p.torsion.for.hyperbolic.satellites}
Let $g = (f_1,\cdots,f_n) \splice L$ where $n \geq 1$, and $L$ a
hyperbolic KGL. If for all $j \in \{1,2,\cdots,n\}$,
$H_i\K_{3,1}(f_j)$ contains no elements of order $p$ for all $i <
2p-2$, then $H_i \K_{3,1}(g)$ contains no elements of order $p$ for
$i < 2p-2$.
\end{lem}

\begin{proof}
In this case, there is a homotopy equivalence $$\K_{3,1}(g) \simeq
S^1 \times \left( SO(2) \times_{A_g} \prod_{i=1}^n \K_{3,1}(f_i)
\right)$$ where $A_g$ is a cyclic group acting via permutations on
the factors in $\prod_{i=1}^n \K_{3,1}(f_i)$.

To determine whether there is $p$-torsion in the homology of
$\K_{3,1}(g)$,
it suffices to determine the $p$-torsion in case $A_g$ is replaced
by the $p$-Sylow subgroup of $A_g$ given by $H= \cyc{p^n}$ as the
induced map $$S^1 \times \left( SO(2) \times_{H} \prod_{i=1}^n
\K_{3,1}(f_i) \right) \to S^1 \times \left( SO(2) \times_{A_g}
\prod_{i=1}^n \K_{3,1}(f_i) \right)$$ induces a split epimorphism on
the $p-torsion$ subgroup by a classical transfer argument ala'
Cartan-Eilenberg.

Consider the covering map
$$\left( SO(2) \times_{H} \prod_{i=1}^n \K_{3,1}(f_i) \right)
\to \left( SO(2) \times_{A_f} \prod_{i=1}^n \K_{3,1}(f_i) \right)$$
for which the group of covering translations is abelian with group
of covering transformations $A_f / H$.

Homologically, this map is onto the $p$-torsion elements of $H_*
\left( SO(2) \times_{A_f} \prod_{i=1}^n \K_{3,1}(f_i) \right)$ since
the composite of the transfer map with the covering map:
$$H_* \left( SO(2) \times_{A_f} \prod_{i=1}^n \K_{3,1}(f_i) \right)
\to H_* \left( SO(2) \times_{A_f} \prod_{i=1}^n \K_{3,1}(f_i)
\right)$$ is multiplication by $|\frac{A_f}{H}|$ which is coprime to
$p$.

The Lemma follows at once from Lemma \ref{lem:bundles.over.circles}.
\end{proof}

\section{On the subspace generated by cabling and summation}
\label{tcns}

The purpose of this section is to describe the subspace $\mathcal
T\K_{3,1}$ of $\K_{3,1}$, consisting of the path components of
$\K_{3,1}$ containing the unknot and all knots generated from the
unknot by iterating the cabling and connected-sum operations. An
alternate description of the space $\mathcal T\K_{3,1}$ is that it
consists of precisely those long knots whose complements have
JSJ-decompositions containing only Seifert-fibred manifolds.

First define the space $$\mathcal T = \amalg_{1<p<q, (p,q) = 1}
\K_{3,1}(f(p,q))$$ where $f(p,q)$ denotes a $(p,q)$-torus knot. Thus
$\K_{3,1}(f(p,q))$ has the homotopy type of $S^1$. Consider the
James construction
$$J(\mathcal T \amalg \{*\}) = \amalg_{0 \leq n} {\mathcal T}^n$$ with
$\mathcal T ^0 = \{*\}$, the base-point. Write
$$J_{\K} = \amalg_{1 \leq n} {\mathcal T}^n$$ with
$$J(\mathcal T \amalg \{*\}) = J_{\K} \amalg \{*\}.$$

Spaces $Y_n$ are specified inductively in terms of $J_{\K}$ as
follows.
\begin{itemize}
  \item[*] $Y_0 = \Cu_2(J_{\K}\amalg \{*\})$ and
  \item[*] $Y_{n+1} = \{([\Cu_2(Y_n)]-Y_n) \times J_{\K})\} \amalg
  Y_n$.
\end{itemize} Notice that $Y_{n}$ is naturally a subspace
of $Y_{n+1}$ and that all of these may be regarded as subspaces of
$\K_{3,1}$ in the following way: There are induced maps
$$\K_{3,1}\times J(\mathcal T \amalg \{*\}) \to \K_{3,1}$$ induced
by cabling.

Define $$ \mathcal T\K_{3,1} = \cup_{n \geq 0} Y_n.$$ Notice that
$\mathcal T\K_{3,1} $ is the subspace of $\K_{3,1}$ which contains
the path-components of $p/q$-torus knots and which is closed under
the operations of cabling and sums. That is, there are induced maps
$\K_{3,1}\times J(\mathcal T \amalg \{*\}) \to \K_{3,1} $ induced by
cabling. There is an induced inclusion $\mathcal T\K_{3,1}  \to
\K_{3,1}$.

The remainder of this section gives features of the homology of
$\mathcal T\K_{3,1} $. First notice that homology commutes with
inductive co-limits and so there are isomorphisms
$$\varinjlim H_*(Y_n) \to\ H_*(\varinjlim Y_n)\to\ H_*(\mathcal T\K_{3,1} ).$$
Recall the construction $\Gamma(X)$ as given in section
\ref{c2homol}.

To describe the homology of $Y_n$, restrict to field coefficients.
$\F$, Recall the natural splitting of graded vectors spaces
$$\bar H_*(X; \mathbb F) \oplus \Gamma(X; \mathbb F) \to\
 H_*(\Cu_2(X); \mathbb F))$$ for a choice
of graded vector space $$\Gamma(X; \mathbb F)$$ which is functor of
$H_*(X; \F)$ as listed in section \ref{c2homol}. Thus $H_*(Y_{n+1})$
is given in terms of the construction $\Gamma(X; \mathbb F)$ in case
$X = Y_n$.

\begin{prop}\label{subspace generated by torus knots}
The natural map $H_*(Y_n) \to\ H_*(Y_{n+1})$ is a split
monomorphism.
\end{prop}

Notice that the homology of $\mathcal T\K_{3,1}$ exhibits a
fractal-like behaviour reflecting the geometry in Budney's theorem
\cite{B} and iterations of the constructions $\Gamma(X)$ as given in
section \ref{c2homol}. Namely, this homological behaviour arises by
first considering $Y_1 = \Cu_2(J_{\K})$ together with the homology
$H_*\left( \Cu_2(X \amalg \{*\});\F\right)$ as given in section
\ref{c2homol} as follows where
$$V = H_*(X;\F).$$
\begin{enumerate}
    \item $S[\sigma^{-1}L[\sigma(V]]$ if the characteristic of $\F$ is $0$,
    \item $S[\sigma^{-1}L^{(2)}[\sigma(V)]]$ if $\F = \F_2$ and
    \item $S[\sigma^{-1}L^{(p)}[\sigma(V)] \oplus \sigma^{-2}W^p[\sigma(V)]]$
    if $\F = \F_p$ for odd primes $p$.
\end{enumerate}

On a simpler note, let $\aleph \K_{3,1}$ denote the subspace of $\K_{3,1}$
consisting of: all unknots, torus knots, and all connect-sums of torus knots.
Thus, $\aleph \K_{3,1}$ is a $2$-cubes subspace of $\K_{3,1}$ and
$$\aleph \K_{3,1} \simeq \Cu_2\left(\{*\} \sqcup \bigsqcup_{\Zed} S^1\right).$$
By May \cite{M}, 
$$B(\aleph \K_{3,1}) \simeq \Omega^2 \Sigma^2\left(\{*\} \sqcup \bigsqcup_{\Zed} S^1\right)$$
which has the homotopy-type of
$$\Omega^2 \left( \bigvee_\Zed \left( S^2 \vee S^3 \right)\right)$$
where the union and wedge index set $\Zed$ corresponds to the isotopy classes of torus knots.
Thus, by the Hilton-Milnor theorem the homotopy groups of $B (\aleph \K_{3,1})$ contain the homotopy
groups of all spheres (of dimension $\geq 2$) in profusion.

\section{Closed Knots and Homology}\label{closedhom}

The purpose of this section is to use results of the earlier
sections to give information about the space $\Emb(S^1,S^3)$. Recall
the homeomorphism $$S^3 \times \Emb_*(S^1,S^3)\to \Emb(S^1,S^3)$$ for
which $Emb_*(S^1,S^3)$ denotes the pointed embeddings. Thus there
are isomorphisms $$H_*(\Emb(S^1,S^3)) \to H_*(S^3) \otimes
H_*(\Emb_*(S^1,S^3))$$ by Proposition \ref{neq3long_closed_rel}.

Information giving the structure of the bundle $\K_{3,1} \to
\Emb_*(S^1,S^3) \to S^2$ was worked out earlier. That structure is
used to provide information concerning $H_*(\Emb_*(S^1,S^3))$ by a
Mayer-Vietoris argument.

Let $D_1$ and $D_2$ be two discs in $S^2$ whose union is $S^2$ and
whose intersection is $S^1$. Let $A_1$ and $A_2$ be the preimages of
$D_1$ and $D_2$ under the projection map $\Emb_*(S^1,S^3) \to S^2$,
then both $A_1$ and $A_2$ are homeomorphic to $\K_{3,1} \times D^2$.
Consider the Mayer-Vietoris sequence for $\Emb_*(S^1,S^3) = A_1 \cup
A_2$ where $A_1 \cap A_2$ is homeomorphic to $S^1 \times \K_{3,1}$.

Identify $A_1 \equiv D^2 \times \K_{3,1}$ and $A_2 \equiv D^2 \times
\K_{3,1}$ then the gluing map from $\partial A_1 = S^1 \times
\K_{3,1} \to \partial A_2 = S^1 \times \K_{3,1}$ is the map $S^1
\times \K_{3,1} \ni (t,x) \longmapsto (t,t^2.x) \in S^1 \times
\K_{3,1}$. Thus the Meyer-Vietoris sequence has the form:
$$\cdots \to H_*(S^1 \times \K_{3,1}) \to H_* (D^2\times \K_{3,1})\oplus H_*(D^2\times \K_{3,1}) \to H_*(\Emb_*(S^1,S^3)) \to \cdots$$

where $H_*(S^1 \times \K_{3,1})$ is identified with $H_* \K_{3,1}
\oplus H_{*-1} \K_{3,1}$ and $H_*(D^2 \times \K_{3,1})$ is
identified with $H_*(\K_{3,1})$. The map $H_n \K_{3,1} \oplus
H_{n-1} \K_{3,1} \to H_n(\K_{3,1}) \oplus H_n(\K_{3,1})$ is given by
the $2\times 2$ matrix
$\left( \begin{array}{cc} I & 0 \\
I & 2\kappa_n \end{array}\right)$ where $\mu : SO(2) \times \K_{3,1}
\to \K_{3,1}$ is the $SO(2)$-action on $\K_{3,1}$ and $\kappa_n :
H_{n-1} \K_{3,1} \to H_n \K_{3,1}$ satisfies
$\kappa_n(x)=\mu_*(SO(2) \times x)$

\begin{cor}\label{sesfors1}
There is a short exact sequence
$$ 0 \to coker(2 \kappa_n) \to H_n \Emb_*(S^1,S^3) \to ker(2 \kappa_{n-1}) \to 0.$$
\end{cor}

\begin{cor}
A knot $f : S^1 \to S^3$ is the unknot if and only if its component
in $\Emb(S^1,S^3)$ contains no torsion in its homology. Moreover,
the component of a non-trivial knot in $\Emb(S^1,S^3)$ always has
$2$-torsion in its integral homology.
\begin{proof}
If $f$ is the unknot, the component of $f$ has the homotopy type of
$V_{4,2} = S^3 \times S^2$ which has no torsion in its homology.

If $f$ is non-trivial, then first consider its long knot component,
$\K_{3,1}(f)$. This has a $\Zed$ embedded in its its fundamental
group, embedded as the $2\pi$ rotation around the long axis
\cite{GramainPi1}. We call the embedding $\Zed \to \pi_1
\K_{3,1}(f)$ the Gramain map. In \cite{B2} it's shown that there is
a map $\K_{3,1}(f) \to S^1$ which when composed with the Gramain map
is not null-homotopic.
It follows that $H_1 \K_{3,1}(f)$ contains a copy of the integers,
generated by the Gramain element.  Thus since the image of
$2\kappa_1$ is generated by twice the Gramain element,
$coker(2\kappa_1)$ must contain $2$-torsion.
\end{proof}
\end{cor}

The short exact sequence in Corollary \ref{sesfors1} is not ideal
because it leaves us with extension problems.  We show how the
extension problems can be solved using the techniques of Section
\ref{htpe}.

Observe that, if $f \in \K_{3,1}$ is a prime knot, then there is an
$SO(2)$-equivariant homotopy-equivalence $\K_{3,1}(f) \simeq SO(2)
\times X(f_i)$ where the $SO(2)$-action on $SO(2) \times X(f)$ is a
product action, given by left-multiplication on $SO(2)$ and the
trivial action on $X(f)$ (here $X(f)$ is just $\K_{3,1}(f)/SO(2)$).
So prime knot components of $Emb(S^1,S^3)$ have the homotopy type of
$S^3 \times SO(3) \times X(f)$.  As mentioned earlier, the unknot
component has the homotopy-type of $S^3 \times S^2$.

We now investigate the case of a connected-sum of $n \geq 2$ prime
knots, $f = f_1 \# \cdots \# f_n$.  By the above argument, we can
assume $\K_{3,1}(f_i) \simeq SO(2) \times X(f_i)$ for $X(f_i) =
\K_{3,1}(f_i)/SO(2)$. Thus, the component corresponding to $f$ in 
$\Emb_*(S^1,S^3)$ has the homotopy-type of $\Cu_2(n)
\times_{\Sigma_f} \left( \left( SO(3) \times_{SO(2)} SO(2)^n \right)
\times \prod_{i=1}^n X(f_i) \right)$

We determine the homotopy-type of $SO(3) \times_{SO(2)} SO(2)^n$ as
a $\Sigma_n$-space. Consider $SO(2)^n$ to be $\Real^n / \Zed^n$.  Let $D
\subset \Real^n$ be the diagonal $D=\{(t,t,\cdots,t) : t \in \Real
\}$.  Let $P \subset \Real^n$ be the perp of $D$, ie: $P = \{
(x_1,x_2,\cdots,x_n) : \sum_{i=1}^n x_i = 0 \}$. Thus $P / (P \cap
\Zed^n)$ is an $(n-1)$-dimensional torus, which we will denote
$P(n)$. We also define a subgroup $Z(n) \subset P(n)$. The integer
lattice $\Zed^n$ projects (orthogonally) onto a subgroup of $P$, we
further take the image of this subgroup under the quotient map $P
\to P(n)$ and denote this image $Z(n)$. There is a naturally defined
homomorphism $Z(n) \to SO(2)$ given by considering the embedding
$P(n) \to SO(2)^n \equiv \Real^n / \Zed^n$.  For every element $z
\in Z(n)$ there is a unique element $t \in SO(2)$ so that $t.z =0
\in SO(2)^n$.

\begin{prop} Provided $n \geq 1$,
$$SO(3) \times_{SO(2)} SO(2)^n \simeq SO(3) \times_{Z(n)} P(n)$$
where the action of $Z(n)$ on $SO(3)$ is given by the homomorphism
$Z(n) \to SO(2)$. Moreover, this is an $\Sigma_n$-equivariant
homeomorphism where the action of $\Sigma_n$ on $SO(3) \times_{Z(n)}
P(n)$ is a product action, trivial on $SO(3)$ and the natural action
on $P(n) \subset SO(2)^n$.  Since the homomorphism $Z(n) \to SO(2)$
is null-homotopic (as a continuous function), the above bundle is
homeomorphic to a product $SO(3) \times_{Z(n)} P(n) \simeq SO(3)
\times (P(n)/Z(n))$.
\end{prop}

\begin{cor} If $f$ is a connected-sum of $n$ prime knots $n \geq 2$, then
$$SO(3) \times_{SO(2)} \K_{3,1}(f) \simeq SO(3) \times \Cu_2(n) \times_{\Sigma_f}
\left( P(n)/Z(n) \times \prod_{i=1}^n X(f_i) \right)$$ where
$\K_{3,1}(f_i) = SO(2) \times X(f_i)$, so the fibration $SO(3)
\times_{SO(2)} \K_{3,1}(f) \to S^2$ is just projection onto $SO(3)$
then onto $S^2$.
\end{cor}

We now perform the analogous computations for $\Emb(S^1,\Real^3)$. 
Proposition \ref{long-cl} gives us the analogous bundle 
$C \rtimes \K_{3,1} \to \Emb(S^1,\Real^3) \to S^2$.
Decomposing $S^2$ as the union of two discs, one gets a
Meyer-Vietoris sequence

$$\cdots \to H_*(S^1 \times (C \rtimes \K_{3,1})) \to H_* (C \rtimes \K_{3,1})\oplus 
H_*(C \rtimes \K_{3,1}) \to H_*(\Emb(S^1,\Real^3)) \to \cdots$$

which splits into short exact sequences as in Corollary \ref{sesfors1}:

$$ 0 \to coker(2\kappa_n') \to H_n\Emb(S^1,\Real^3) \to ker(2\kappa_{n-1}') \to 0$$
where $\kappa_n' : H_{n-1}(C \rtimes \K_{3,1}) \to H_n(C \rtimes \K_{3,1})$ is
given by $\mu_*(SO(2) \times \cdot)$ for the $SO(2)$ action $\mu$ on $C \rtimes \K_{3,1}$.
The bundle $C_f \to C_f \rtimes \K_{3,1}(f) \to \K_{3,1}$ is split, and the monodromy
acts trivially on $H_* C_f$ since $C_f$ is a homology $S^1$ with $H_1(C_f)$ generated
by a meridional curve. Thus, $H_* (C \rtimes \K_{3,1}) \simeq H_* S^1 \otimes H_* \K_{3,1}$ and
$H_n(C\rtimes \K_{3,1}) = (H_0 S^1 \otimes H_n\K_{3,1}) \oplus (H_1 S^1 \otimes H_{n-1}\K_{3,1})$.

$\kappa_n'$ has a description in terms of $\kappa_n$ and $\kappa_{n-1}$.  
Let $\alpha_i \in H_i S^1$ represent
the standard generators of $H_i S^1$ for $i=0,1$.  Then 
$\kappa_n'(\alpha_0 \otimes x) =  \alpha_1 \otimes x + \alpha_0 \otimes \kappa_n(x)$ and
$\kappa_n'(\alpha_1 \otimes x) = -\alpha_1 \otimes \kappa_{n-1}(x)$. Thus, 
$\kappa_n'$ can be thought of as a map $\kappa_n' : H_{n-2} \K_{3,1} \oplus H_{n-1} \K_{3,1}
\to H_{n-1} \K_{3,1} \oplus H_{n} \K_{3,1}$  given by 
$\kappa_n'(x,y) = (-\kappa_{n-1}(x)+y, \kappa_n(y))$.
Since $\kappa_n \circ \kappa_{n-1} = 0$, $ker(2 \kappa_n')$ is given by the solutions to 
the equation $-2\kappa_{n-1}(x)+2y=0$ for $(x,y) \in H_{n-2} \K_{3,1} \oplus H_{n-1} \K_{3,1}$. 
Thus, 
$$ker(2\kappa_n') \simeq H_{n-2} \K_{3,1} \oplus \tau_2 H_{n-1} \K_{3,1}$$
where if $A$ is an abelian group and $p$ an integer, $\tau_p A$ is the subgroup of $A$ killed
by multiplication by $p$. Similarly,
$$coker(2\kappa_n') \simeq H_{n-1} \K_{3,1} / 2H_{n-1} \K_{3,1} \oplus H_n \K_{3,1}$$

\begin{prop} There is a short exact sequence
$$ 0 \to H_{n-1} \K_{3,1} / 2H_{n-1} \K_{3,1} \oplus H_n \K_{3,1} \to H_n\Emb(S^1,\Real^3) \to
H_{n-3} \K_{3,1} \oplus \tau_2 H_{n-2} \K_{3,1} \to 0 $$
Thus, the component of the unknot in $\Emb(S^1,\Real^3)$ is the unique component
such that its first homology group is torsion.  It is also the unique component so that
its 2nd homology group is trivial.
\begin{proof}
The short exact sequence follows from the above observations.   

That $H_1$ of a non-trivial component
is non-torsion follows from Proposition 6.1 of \cite{B2} and the above short exact sequence.
That $H_1$ of the unknot component is torsion follows from Proposition \ref{long-cl}
and the results in Section \ref{htpe}.  Thus, the component of the unknot in $\Emb(S^1,\Real^3)$
has the homotopy-type of $SO(3)$, and $H_1 SO(3) \simeq \Zed_2$.

The statement about $H_2$ follows from Proposition 6.1 of \cite{B2} 
and the above short exact sequence. 
\end{proof}
\end{prop}

Proposition \ref{prop:Two torsion in h1} gives a criteria for testing 
whether a long knot f is the `long unknot' which can be modified to 
to compare two arbitrary long knots. That procedure involves forming a 
`difference' to be defined precisely. Namely, it is not the case that an 
arbitrary long knot admits an inverse. It is necessary to pass to a setting for
which the inverses exists in order to `take differences'. That
`world' is the group completion of $\K_{3,1}$,
$$\Omega(B\K_{3,1}).$$ This last space admits inverses up to
homotopy. Given two long knots $f$ and $g$, consider their classes in
$\Omega(B\K_{3,1})$ denoted $[f]$, and $[g]$ respectively. Next,
consider the product $$[f] \cdot [g]^{-1} \in \Omega(B\K_{3,1}).$$
The path-component of $[f] \cdot [g]^{-1} \in \Omega(B\K_{3,1})$ has
vanishing first homology group if and only if  $f$ and $g$ are in
the same path-component of $\K_{3,1}$.

\section{Problems}\label{Problems}

The purpose of this section is to list problems which arise
naturally from the work above.

(1): Interpret the rational cohomology of $\K_{3,1}$ in terms of
iterated integrals in the sense of Kohno-Kontsevich-Chen.

(2): Compare the Vassiliev invariants of braids as studied by
T.~Kohno \cite{K2} and the Lie algebra obtained from the descending
central series for the fundamental groups of the spaces $\Cu_2(n)
\times \K_{3,1}({f_1}) \times \cdots \times \K_{3,1}({f_n})$ as well
as the induced invariants for $\Cu_2(n) \times_{\Sigma_f} \left(
\K_{3,1}({f_1}) \times \cdots \times \K_{3,1}({f_n}) \right) $.

(3): A natural connection between the space of long knots and the
mod-$2$ Steenrod algebra arises from the ``group completion'' of
$\K_{3,1}$ given $\Omega B(\K_{3,1})$ \cite{M}. Notice that the
collapse map $X_{\K}\to \{*\}$ induces a map $p:C(\Real^2,X_{\K}
\amalg \{*\} ) \to C(\Real^2,S^0 )$ and there is an induced map
$p:C(\Real^2,X_{\K} \amalg \{*\} ) \to \Zed \times BO$ induced by
the regular representation bundle. Thus there are maps

$$ \Omega B(\K_{3,1}) \to \Omega^2_0 S^2 \to BO $$
with composite denoted $\phi:\Omega B(\K_{3,1}) \to \Zed \times BO$.
The Thom spectrum of $\phi$ is a wedge of Eilenberg-Mac Lane spectra
$H\F_2$ and thus the mod-$2$ co-homology of the Thom spectrum
$M\Omega B(\K_{3,1})$ is free over the mod-$2$ Steenrod algebra.
Interpret the Steenrod operations in terms of knots.

(4): The Goodwillie Calculus mapping space models $AM_j(\I^n)$ for
$\K_{n,1}$ constructed by Sinha \cite{Sinha} have a natural
homotopy-associative pairing. This pairing makes $\pi_0 AM_3(\I^3)$
into a group, isomorphic to the integers. As an invariant of knots,
$\pi_0 \K_{3,1} \to \pi_0 AM_3(\I^3) \simeq \Zed$ is the essentially
unique type-2 finite-type invariant of knots \cite{B4}. This raises
the question, do the maps $\K_{n,1} \to AM_j(\I^n)$ factor through
the group completion, $\K_{n,1} \to \Omega B\K_{n,1} \to
AM_j(\I^n)$?

(5): Combine the structures here with that of Khovanov homology.

\end{document}